# Forecasting Outside the Box: Application-Driven Optimal Pointwise Forecasts for a Class of Two-Stage Stochastic Programs[*]

Tito Homem-de-Mello[†1], Juan Valencia[2], Felipe Lagos[3], and Guido Lagos[4]

[1]*School of Business, Universidad Adolfo Ibáñez, Chile*
[2]*School of Business, Universidad Adolfo Ibáñez, Chile*
[3]*Faculty of Engineering and Sciences, Universidad Adolfo Ibáñez, Chile*
[4]*Faculty of Engineering and Sciences, Universidad Adolfo Ibáñez, Chile*

**Abstract**

We study a class of two-stage stochastic programs, namely, those with fixed recourse matrix and fixed costs, and linear second stage. We show that, under mild assumptions, the problem can be solved with just *one* scenario, which we call an "optimal scenario." Such a scenario does not have to be unique and may fall outside the support of the underlying distribution. Although finding an optimal scenario in general might be hard, we show that the result can be particularly useful in the case of stochastic optimization problems with *contextual information*, where the goal is to optimize the expected value of a certain function given some contextual information (e.g., previous demand, customer type, etc.) that accompany the main data of interest. The contextual information allows for a better estimation of the quantity of interest via machine learning methods. We focus on a class of learning methods—sometimes called in the literature decision-focused learning—that integrate the learning and optimization procedures by means of a bilevel optimization formulation, which determines the parameters for pointwise forecasts. By using the optimal scenario result, we prove that when such models are applied to the class of contextual two-stage problems considered in this paper, the pointwise forecasts computed from the bilevel optimization formulation actually yield asymptotically the best approximation of an optimal scenario within the modeler's pre-specified set of parameterized forecast functions. Numerical results conducted with inventory problems from the literature (with synthetic data) as well as a bike-sharing problem with real data demonstrate that the proposed approach performs well when compared to benchmark methods from the literature.

## 1 Introduction

The area of stochastic optimization has evolved considerably in the past decade. Traditionally, stochastic optimization models assumed the existence of a known probability distribution to represent the underlying uncertainty, and formulated the problem in terms of optimizing the expected value (or another risk measure) of a certain function with respect to a decision variable, where the expectation corresponds to that distribution.

---

One can formulate the problem in a generic way as

$$\min_{z \in Z} \mathbb{E}_P\left[G(z, \xi)\right], \tag{1}$$

where $z$ is the decision variable, $\xi$ represents the uncertainty, and $P$ is the distribution of $\xi$. Much of the earlier efforts aimed at developing *scenario generation/reduction techniques* for the case where $P$ has either large or infinite support, in order to approximate the original problem with one in which the distribution has moderately-sized support and hence decomposition methods can be used to solve the problem. Among the scenario generation methods are clustering (Dupačová, Consigli, and Wallace, 2000), moment-matching techniques (Hoyland and Wallace, 2001; Hoyland, Kaut, and Wallace, 2003; Mehrotra and Papp, 2014), and Monte Carlo/Quasi-Monte Carlo methods (Homem-de-Mello and Bayraksan, 2014; Leövey and Römisch, 2015; Shapiro, Dentcheva, and Ruszczynski, 2021, Chapter 5). Another class of methods, based on probability metrics, aims at finding a distribution $Q$ with relatively few scenarios in such a way that $Q$ minimizes a distance $d(P, Q)$ between $Q$ and the original distribution $P$. We refer to Dupačová, Gröwe-Kuska, and Römisch (2003); Heitsch and Römisch (2003); Heitsch and Römisch (2009); Pflug (2001); Pflug and Pichler (2011) and references therein for further discussions on this type of methods.

Recent efforts have been directed to developing scenario generation/reduction techniques that use the information from the optimization problem at hand rather using just the distribution of the underlying random variables. Such methods typically focus on two-stage stochastic programs in order to exploit the structure of the problem. For instance, Bertsimas and Mundru (2023) define the distance between two probability distributions in terms of the cost functions of the optimization problem. Henrion and Römisch (2022) use problem information to compute the distribution that leads to the best uniform approximation of $\mathbb{E}_P\left[G(z, \xi)\right]$ over all feasible $z$, whereas Zhang, Wang, Jacquillat, and Wang (2023) propose a scenario subset selection model that optimizes the approximation of the recourse function over a pool of first-stage solutions. Keutchayan, Ortmann, and Rei (2023) develop a problem-driven scenario clustering method that produces a partition of the scenario set that enables representative scenarios to be identified. A different approach is used in Arpón, Homem-de-Mello, and Pagnoncelli (2018) and Fairbrother, Turner, and Wallace (2019), who work with a variation of (1) where the objective function is a tail risk measure such as Conditional Value-at-Risk, and develop scenario generation methods that exploit the structure of that objective.

In this paper we consider problems of the form (1) when the function $G$ corresponds to a *two-stage stochastic program*. More specifically, we consider two-stage stochastic programs with fixed recourse and fixed costs (henceforth denoted FRFC), i.e.,

$$\min_{z \in Z} f(z) + \mathbb{E}\left[Q(z, \xi)\right] \tag{2}$$

where $Z \subseteq \mathbb{R}^n$ is a convex set, $f : \mathbb{R}^n \mapsto \mathbb{R}$ is a convex function such that $Z \subseteq \text{Dom } f$, $Q$ is the second stage function

$$\begin{aligned}
Q(z, \xi) = \quad & \min \; q^\top y & (3) \\
\text{s.t.} \quad & Wy = h - Tz & (4) \\
& y \geq 0 & (5)
\end{aligned}$$

and $\xi$ denotes the random element $(h, T)$. The qualifier "fixed recourse and fixed costs" refers to the fact that in the second-stage problem (3)-(5) neither the matrix $W$ nor the vector $q$ are random. The class of two-stage stochastic programs with FRFC is actually quite large, as it includes problems where the uncertainty corresponds to demand, as is the case in many inventory, energy, capacity planning and logistics problems; in fact, the vast majority of two-stage test problems reported in the literature, or available in public repositories such as `stoprog.org`, are of FRFC type. A notable exception consists of problems where the coefficient $q$ corresponds to random prices or returns, as in the case of portfolio models for example.

A key result of the paper is the proof that, under mild assumptions, *the two-stage stochastic program given in* (2)-(5) *can be solved with only one scenario*, in the sense that there exists a scenario $\xi^* = (h^*, T^*)$ (possibly outside the support of $\xi$)—which we call an "optimal scenario"—such that solving the problem $\min_{z \in Z} f(z) + Q(z, \xi^*)$ yields an optimal solution that is also optimal for (2)-(5). In other words, it suffices to solve the simpler problem

$$\min_{z \in Z} \ f(z) + q^\top y \tag{6}$$

$$\text{s.t.} \quad T^* z + W y = h^* \tag{7}$$

$$y \geq 0. \tag{8}$$

instead of (2)-(5). This is a surprising result, which to best of our knowledge has not been shown in the literature. In fact, one of the main arguments for solving stochastic optimization problems such as (1) instead of the simpler one-scenario problem $\min_{u \in U} \ G(u, \bar{\xi})$ for some fixed scenario $\bar{\xi}$ for some fixed scenario $\xi^*$ is the fact that the one-scenario problem does not capture the variability of the uncertainty. For instance, Wallace (2000) presents a compelling argument by means of very simple examples of stochastic optimization problems for which solving the problem for one scenario (regardless of the choice of the scenario) can never yield the same solution as the stochastic one.

So, how to reconcile our main result with the proven need to use the full distribution of the uncertainty? As we shall see later in the paper, the key lies in the particular characteristics of two-stage stochastic programs with FRFC. In light of that result, it is not surprising that none of examples presented in Wallace (2000) are two-stage stochastic programs with FRFC, so there is no contradiction between our result and the conclusions in Wallace (2000)[1]. Interestingly, our result can also be viewed as a generalization of a property that holds for the well-known newsvendor problem, which is a particular case of a two-stage stochastic program with FRFC, as will be discussed in Section 2.

The one-scenario result, while appealing from a theoretical perspective, has its caveats. One is that an optimal scenario may not respect dependencies (either functional or statistical ones) among the random variables; we will discuss this issue in more detail later. Another caveat is that it is a result about *existence* of an optimal scenario; finding one such scenario may be difficult. On the other hand, we make no claims about uniqueness of the optimal scenario. For our purposes, it suffices to know that we can search for *some* scenario with the property that solving the corresponding one-scenario problem—which is just a deterministic problem, with no random variables— yields the same solution as the stochastic one.

Those caveats notwithstanding, the result might be useful in a few ways. For example, it is conceivable developing an algorithm that searches over the space of scenarios instead of over the space of decisions, which could be advantageous in some settings; we shall see an example of that in Section 7. Also, optimal scenarios may have an interesting interpretation for the decision maker, as we will illustrate with an example in Section 2.

We also exploit the consequences of the one-scenario result in the setting where there is *contextual information*. The enormous growth in the availability of data in recent years, and more specifically the presence of contextual information in the data—also called *covariates*, or *features* in the literature— has led to the development of new models in stochastic optimization. In such models the uncertainty represented by $\xi$ can be predicted to some extent by the available contextual information. Thus, the goal is to optimize the expected value of a certain function conditionally on some given value of the contextual information, henceforth called the *contextual information of interest.* That is, the goal is to find the best decision corresponding to data with some characteristics. An example of such a situation will be discussed in Section 7; the goal is there to determine the best assignment of bicycles on a given day for a bike-sharing service, given some contextual information such as the weather forecast for that day. Formally, we can write the problem as

$$\min_{z \in Z} \mathbb{E}\left[G(z, \xi) \mid X = x\right], \tag{9}$$

where again $z \in \mathbb{R}^n$ is the decision variable, $Z \subseteq \mathbb{R}^n$ is the feasibility set, $\xi$ is a random variable with support $\Xi \subseteq \mathbb{R}^m$ that represents the main uncertainty, $X$ is a random variable with support $\mathcal{X} \subseteq \mathbb{R}^s$ that represents the contextual information and $x$ is the contextual information of interest.

---
[1] Actually, Example 1 in Wallace (2000) could be formulated as a very special case of a two-stage stochastic program with FRFC where, in the notation of (3)-(5), the function $Q$ is either 0 or $\infty$. Such a model however does not have relatively complete recourse, a common assumption that we also make in our developments.

The difficulty of optimizing by taking into account the presence of contextual information as in (9) is how to define a proper *conditional* distribution of $\xi$ given $X = x$ to use in the model. Simple techniques such as "slicing" the data to keep only the data points corresponding to the contextual information of interest are not practical since the resulting dataset may be too small or even empty, which occurs in case the particular contextual information of interest has not been observed in the dataset. An alternative in such cases is to use *machine learning* approaches that can forecast the uncertainty as a function of the contextual information. Modern machine learning methods such as neural networks, regularized regression and classification trees, among others, can be used to learn the dependence of $\xi$ on $x$ (see, e.g., Bottou, Curtis, and Nocedal (2018); Hastie, Tibshirani, and Friedman (2009)). The key then is how to combine these predictions with the optimization model.

The use of machine-learning-based forecasts for optimization can be accomplished in multiple ways, and has become a fruitful topic in stochastic optimization with many papers in the past five years. As discussed in the recent survey paper by Sadana, Chenreddy, Delage, Forel, Frejinger, and Vidal (2024), three main approaches can be found in the literature: (i) *decision-rule optimization*, which aims to approximate directly the optimal solutions of (9) as a function of $x$ by means of techniques such as linear decision rules or reproducing kernel Hilbert spaces; (ii) *sequential learning and optimization*, which uses machine learning techniques to estimate the conditional distribution of the uncertainty given the contextual information of interest, and then applies standard methods to solve the stochastic optimization problem corresponding to that conditional distribution; and (iii) *integrated learning and optimization*, where the forecast and optimization are combined within the same problem. For instance, in the SPO (smart "predict, then optimize") framework of Elmachtoub and Grigas (2021), for a given feature $x$ the problem is written as

$$\min_{z \in Z} \mathbb{E}\left[\xi^T z \mid X = x\right] \;=\; \min_{z \in Z} (\mathbb{E}\left[\xi|x\right])^T z, \tag{10}$$

where we use $\mathbb{E}\left[\xi|x\right]$ as a short for $\mathbb{E}\left[\xi \mid X = x\right]$. Note that the linearity of the above model implies that in order to solve the problem we only need an estimate $\xi_x$ of $\mathbb{E}\left[\xi|x\right]$, i.e., a pointwise forecast. The key idea of the SPO approach is to measure the *decision error* induced by the estimation error of $\mathbb{E}\left[\xi|x\right]$), and to measure the performance of the prediction in terms of its impact in the objective function instead of using a standard error criterion such as least-squares. Such an idea can actually be traced back to Bengio (1997) but has gained traction in recent years—albeit with different names such as *integrated conditional expectation and optimization* (Grigas, Qi, and Shen, 2021), *end-to-end learning* (Donti, Amos, and Kolter, 2017), *application-* *driven learning* (Dias-Garcia, Street, Homem-de-Mello, and Muñoz, 2024), and *decision-focused learning* (Mandi, Bucarey, Tchomba, and Guns, 2022), in addition to the other terminology mentioned above. The idea has gone even beyond scientific papers; for instance, a recent Harvard Business Review article describes the implementation of a forecasting methodology for supply chains, called *optimal machine learning* by the authors, that "involves using artificial intelligence technology to create a mathematical model that takes data inputs [...] and links them to planning decisions" (Agrawal, Cohen, Deshpande, and Deshpande, 2024).

We study how to use the aforementioned one-scenario result to solve problems of the form (9) when the function $G$ corresponds to a two-stage stochastic program, using the problem information to measure the forecast error. Notice that the problem is considerably harder than (10), since in principle we need, as discussed earlier, to forecast the entire conditional distribution of $\xi$ given $X = x$ using a problem-based approach. This is in fact the approach used by Grigas et al. (2021), although the approach has some limitations such as fixing in advance the support of the distribution. The one-scenario result ensures that forecasting the entire conditional distribution is not necessary for the class of two-stage stochastic programs with FRFC.

Indeed, consider a mapping $z^*(x)$ which yields an optimal solution to (9) as a function of $x$. As determining the entire mapping $z^*(\cdot)$ is impractical from a computational perspective, it is natural to think of ways to *approximate* that function, and this is precisely what we do in this paper by approximating the function that maps $x$ to an optimal scenario for problem (9). That is, instead of approximating the optimal *solution* mapping as done for instance in Ban and Rudin (2019), we approximate an optimal *scenario* mapping, call it $\xi_x^*$ to emphasize its dependence on $x$. Thus, our approach can be viewed as a bridge between the "integrated learning and optimization" and the "decision rule optimization" methods, in the sense that we aim at producing problem-based biased forecasts but we already know that there exists an optimal forecast which is actually a (linear) function of the optimal solution. Note that in our approach we do not need to deal with distributions, only with parameterized pointwise forecasts $\Psi(\theta, x)$ which constitute the vast majority of forecasts obtained with machine learning techniques. Once an approximation to an optimal scenario mapping (call it $\Psi(\theta^*, x)$) is constructed from training data, we can easily obtain the corresponding solution to (9) for any given $x$ by solving the one-scenario problem in (6)-(8) with $\Psi(\theta^*, x)$ in place of $\xi^*$.

Naturally, the task of approximating an optimal scenario mapping is not simple. This is where the notion of *application-driven forecasts* developed in Dias-Garcia et al. (2024) becomes key. As discussed earlier, the approach in that paper falls into the category of works that measure the quality of the pointwise forecast in terms of its impact on the optimization problem; in the case of Dias-Garcia et al. (2024), this is accomplished by solving a bilevel problem (see also Muñoz, Pineda, and Morales 2022 and Morales, Munoz, and Pineda 2023 for similar approaches). A distinctive feature of the approach in Dias-Garcia et al. (2024) that is useful here is the fact that the method aims at finding the *best possible values* of the parameters of the forecast function; thus, as we show in the present paper, as long as the class of forecast functions is flexible enough, the pointwise forecast yielded by the algorithm will be a good approximation of an optimal scenario. There is however a trade-off between the flexibility allowed by the class of forecast functions and the computational effort required to solve the bilevel model; still, such an effort is spent at the *training* stage—as discussed above, once the optimal parameters are found, solving (9) amounts to solving a simple problem. Our numerical results in Section 7, where we study three problems from the literature—two with synthetic data and one with real data—indicate that the one-scenario forecast actually performs very well.

As seen above, prediction methods that rely on context variables and bilevel models have been studied in the literature. Our work, however, introduces the novel result of a single-scenario equivalent formulation, along with a series of new contributions that build upon it. Specifically, we propose a solution strategy in which, for a given context $x \in \mathcal{X}$, a policy yields a solution to the stochastic problem. Our methodology is sufficiently flexible to incorporate nonlinear prediction functions, such as classification and regression trees (CARTs), while still preserving our convergence guarantees. As in Muñoz et al. (2022), Dias-Garcia et al. (2024), and Morales et al. (2023), we also employ KKT conditions to reformulate and simplify the two-level model. To address the resulting problem, we adopt the meta-algorithm proposed in Dias-Garcia et al. (2024). The following sections provide a detailed exposition of these contributions.

The remainder of the paper structured as follows: in Section 2 we present, as a motivating example for the results in the paper, a variation of the classical newsvendor model. In Section 3 we introduce our main result, which is the existence of an optimal scenario for two-stage stochastic programs with FRFC. Section 4 discusses optimal solution mappings as a generalization of (9) for multiple values of the contextual information $x$. In Section 5 we present our proposed approach and show that, in the limit, it produces the best possible parameters for given forecast function $\Psi$. We also show conditions under which the forecast will indeed yield a good approximation of an optimal scenario warranted by the result in Section 3. Section 6 presents two specific algorithms that implement our strategy—one based on regression, the other based on classification and regression trees (CART)—and a few algorithms from the literature that are used as benchmarks. Finally, in Section 7 we show numerical results for both non-contextual and contextual cases, using problems from the literature for the latter. Concluding remarks are presented in Section 8.

## 2 A motivating example: the newsvendor model with unreliable supplier

In this section we illustrate some attributes of the one-scenario approach by means of a simple model for which we can derive analytical solutions. We show that there may be infinitely many optimal scenarios, and that some of these scenarios may fall outside of the original support of the random variables.

Consider the classical newsvendor problem, where a retailer purchases $z$ units of product for a cost of $c$ (per unit) before knowing the random demand $D$, and sells the product for a price $p > c$ per unit. Non-sold items incur an inventory cost of $\eta$ per unit, whereas the penalty per unit for unmet demand is $\pi$. The goal is to find the amount $z$ that minimizes the expected cost. The problem is formulated as

$$\min_{z \geq 0} \mathbb{E}\left[cz - p\min(z, D) + \eta[z - D]_+ + \pi[z - D]_-\right], \tag{11}$$

where $[a]_+ := \max\{a, 0\}$ and $[a]_- := \max\{-a, 0\}$. It is well known that, when the distribution of $D$ is continuous, the optimal solution of the above problem is unique and given by $z^* := F^{-1}(\phi)$, where $F$ is the cumulative distribution function of $D$—assumed to be invertible, for simplicity—and $\phi$ is the so-called *critical ratio*, defined as

$$\phi := \frac{p + \pi - c}{p + \pi + \eta} \tag{12}$$

(see, e.g., Gallego, Ryan, and Simchi-Levi 2001). The expression in (12) can be interpreted as $C_u/(C_u + C_o)$, where $C_u$ is the per-unit underage cost $\pi + (p - c)$ and $C_o$ is the per-unit overage cost $c + \eta$. By noticing that $\min(z, D) = z - [z - D]_+$, it follows that the term inside the expectation in (11) can be written as $(c - p)z + (p + \eta)[z - D]_+ + \pi[z - D]_-$, and so (11) can be written as the two-stage stochastic program

$$\min_{z \geq 0} (c - p)z + \mathbb{E}\left[Q(z, D)\right] \tag{13}$$

where $Q$ is the second stage function

$$Q(z, D) = \min_{y^+, y^-} \ (p + \eta)y^+ + \pi y^- \tag{14}$$

$$\text{s.t.} \quad y^+ - y^- = -D + z \tag{15}$$

$$y^+, y^- \geq 0. \tag{16}$$

In the notation of (3)-(5), we have $h = -D$ and $T = -1$. It is easy to see that, if we solve (11) with a fixed demand $\bar{D}$, the optimal solution of this one-scenario problem is trivially $z = \bar{D}$, so by choosing $\bar{D} = z^* = F^{-1}(\phi)$ we recover the optimal solution of the original problem. That is, for the purposes of solving the problem, the quantile $F^{-1}(\phi)$ actually represents the entire distribution of $D$.

Consider now a variation of the newsvendor problem where the supplier is *unreliable*. That is, given an order quantity $z$, the amount that is actually delivered (and paid for) is $Uz$, where $U$ is a random variable with support in $(0, 1]$. This is one of the cases studied by Dada, Petruzzi, and Schwarz (2007), who actually show that the optimal solution $z$ of the problem is given by the solution of the equation

$$\mathbb{E}\left[UF(Uz)\right] = \phi \mathbb{E}[U], \tag{17}$$

where $\phi$ is the same critical value defined in (12). Note that when the supplier is not unreliable (i.e. $U \equiv 1$), the solution of (17) indeed coincides with that of the standard newsvendor model. In general, (17) does not have an analytical-form solution, but this can be accomplished in particular cases. For instance, suppose that $U$ is Uniform$(0, 1)$ and that $D$ is Uniform$(0, b)$ for some $b > 0$. Then, it is possible to show (after some algebra) that the solution of (17) in that case is

$$z^* := \begin{cases} \frac{3}{2}\phi b & \text{if } \phi \in [0, 2/3] \\ \frac{b}{\sqrt{3(1-\phi)}} & \text{if } \phi \in [2/3, 1) \end{cases} \tag{18}$$

(see Appendix A for the proof). We see that, when $\phi \leq 2/3$, the solution is to order 50% more than in the standard case (which is $F^{-1}(\phi) = \phi b$), up to the maximum demand $b$. When $\phi \geq 2/3$, we have a different expression, and we see that, as $\phi$ approaches 1, $z^*$ goes to infinity. This makes sense—note that $\phi \to 1$ means that the overage cost $C_o$ goes to zero and so the solution is to order as much as possible to compensate for the unreliability of the supplier, which corresponds to a reliability factor $U$ that can be arbitrarily close to zero in this case.

Note that the two-stage formulation of the model with unreliable supplier is similar to that in (13)-(16). However, we need a new second-stage variable to indicate the amount actually delivered by the supplier. We can formulate the problem as

$$\min_{z \geq 0} \mathbb{E}\left[Q(z, D, U)\right] \tag{19}$$

where $Q$ is the second stage function

$$Q(z, D, U) = \min_{y^+, y^-, v} \ (c - p)v + (p + \eta)y^+ + \pi y^- \tag{20}$$

$$\text{s.t.} \quad y^+ - y^- - v = -D \tag{21}$$

$$v = Uz \tag{22}$$

$$y^+, y^-, v \geq 0. \tag{23}$$

Similarly as before, we can see that, if we solve this problem with a fixed demand $\bar{D}$ and a fixed reliability level $\bar{U} > 0$, the optimal solution of this one-scenario problem is the one that makes $y^+ = y^- = 0$, i.e. $z = \bar{D}/\bar{U}$. Thus, by fixing any $\bar{U} > 0$ and choosing $\bar{D} = \bar{U}z^*$—with $z^*$ defined in (18)—we recover the optimal solution of the original problem, i.e., the pair $(\bar{D}, \bar{U})$ chosen in this fashion is an optimal scenario.

The case of unreliable supplier offers some new insights about optimal scenarios. First, as seen above, any one-scenario problem with a pair $(\bar{D}, \bar{U})$ such that $\bar{U} > 0$ and $\bar{D} = \bar{U}z^*$ yields the same solution as the original problem; thus, there may be infinitely many optimal scenarios. Also, the condition $\bar{D} = \bar{U}z^*$ shows that, if one fixes the reliability variable $U$ to some value (for example, its mean), then for the purposes of solving the problem, the value given by $\bar{U}z^*$ actually represents the entire distribution of $D$. However, such a value can actually be *outside* the support of the distribution. For instance, suppose we fix $\bar{U}$ to the mean value $\mathbb{E}[U] = 1/2$, so we have $\bar{D} := z^*/2$. As discussed above, $z^* \to \infty$ as $\phi \to 1$, and in particular we have that $\bar{D} > b$ for $\phi > 11/12$. That is, this optimal scenario overestimates demand even beyond its maximum value in order to compensate for the supplier's unreliability, but the one-scenario problem with $(D, U) \equiv (\bar{D}, \bar{U})$ still yields the same optimal solution $z^*$ as the stochastic problem.

## 3 One-scenario optimality

In this section we state a key result of the paper, which will be used in the sequel. It ensures that, under certain conditions, there exists an optimal scenario for the two-stage problem given in (2)-(5) such that if the problem is solved only with that scenario, it yields an optimal solution to the original problem.

To proceed, we make the following assumptions:

**Assumption 1.** *The feasibility set $Z$ is non-empty, and the function $f$ defined in (2) is such that the relative interior of its domain (denoted $\mathrm{ri}(\mathrm{Dom}\, f)$) is non-empty.*

The assumption on non-emptiness of $Z$ is natural, otherwise the two-stage problem of interest is infeasible. The assumption on $f$ is mild, holding for example if $\mathrm{Dom}\, f = \mathbb{R}^n$.

Consider now the dual problem of (3)-(5):

$$\max\ (h - Tz)^\top u \tag{24}$$
$$\text{s.t.} \quad W^\top u \leq q \tag{25}$$
$$u \in \mathbb{R}^m. \tag{26}$$

**Assumption 2.** *The feasibility set $U$ defined by (25)-(26) is non-empty and bounded.*

This assumption ensures that $Q(z, \xi)$ is finite for all values of $z$ and $\xi$.

We state now the main result. Note that no assumptions are made about the distribution of the uncertainty.

**Theorem 1.** *Suppose Assumptions 1-2 hold, and let $z^*$ be an optimal solution to (2)-(5). Let $\xi^* = (h^*, T^*)$ be defined such that $T^* := \mathbb{E}[T]$ and $h^* := T^*z^*$. Then, we have that*

$$z^* \in \underset{z \in Z}{\mathrm{argmin}}\ f(z) + Q(z, \xi^*). \tag{27}$$

*Thus, if the optimal solution of the* one-scenario *problem $\min_{z \in Z} f(z) + Q(z, \xi^*)$ is unique, then it must coincide with $z^*$.*

*Proof.* Assumption 2 implies that there exist a set $\widetilde{U} := \{u_1, \ldots, u_k\}$ (where each $u_i$ is a vertex of $U$) such that, for any values of $z$ and $\xi$, we have

$$Q(z, \xi) = \max\{(h - Tz)^\top u \,:\, u \in \widetilde{U}\}. \tag{28}$$

It is well known that $Q(\cdot, \xi)$ is convex for all $\xi$. From (28), define the set

$$\widetilde{U}_z(\xi) := \underset{u \in \widetilde{U}}{\mathrm{argmax}}\ (h - Tz)^\top u. \tag{29}$$

It follows (see, e.g., Rockafellar 1970) that the subdifferential set of $Q$ at $(z, \xi)$ w.r.t $z$ is given by

$$\partial_z Q(z, \xi) = -T^\top \mathrm{conv}(\widetilde{U}_z(\xi)) \subseteq -T^\top \mathrm{conv}(\widetilde{U}), \tag{30}$$

and so by applying the expectation operator to both sides of the above relationship we have that

$$\mathbb{E}\left[\partial_z Q(z, \xi)\right] \subseteq -\mathbb{E}[T]^\top \mathrm{conv}(\widetilde{U}). \tag{31}$$

Now, since $Q$ is a convex function, we have that

$$\partial_z \mathbb{E}\left[Q(z, \xi)\right] = \mathbb{E}\left[\partial_z Q(z, \xi)\right] \tag{32}$$

and thus from (31)-(32) we conclude that

$$\partial_z \mathbb{E}\left[Q(z, \xi)\right] \subseteq -\mathbb{E}[T]^\top \mathrm{conv}(\widetilde{U}). \tag{33}$$

Consider now any optimal solution $z^*$ to (2)-(5). Since such a problem is convex, if follows that the optimality condition for $z^*$ is

$$0 \in \partial_z \big(f(z^*) + \mathbb{E}\left[Q(z^*, \xi)\right]\big) + N_Z(z^*) \tag{34}$$

where $N_Z(z)$ denotes the normal cone of $Z$ at $z \in Z$. By Assumption 1 and 2, we have that

$$\emptyset \subseteq \mathrm{ri}(\mathrm{Dom}\ f) \cap \mathrm{ri}(\mathrm{Dom}\ Q(\cdot, \xi)),$$

so the formula $\partial_z \big(f(z) + \mathbb{E}\left[Q(z, \xi)\right]\big) = \partial_z f(z) + \partial_z \mathbb{E}\left[Q(z, \xi)\right]$ applies (Rockafellar, 1970) and then from (33) and (34) we then have that

$$0 \in \partial_z f(z^*) + N_Z(z^*) - \mathbb{E}[T]^\top \mathrm{conv}(\widetilde{U}). \tag{35}$$

Define now $\xi^* = (h^*, T^*)$ such that $T^* := \mathbb{E}[T]$ and $h^* := T^* z^*$. Then, from (29) we have

$$\widetilde{U}_{z^*}(\xi^*) = \widetilde{U} \tag{36}$$

(since any dual solution is optimal in that case) and thus from (30) it follows that

$$\partial_z Q(z^*, \xi^*) = -(T^*)^\top \mathrm{conv}(\widetilde{U}) = -\mathbb{E}[T]^\top \mathrm{conv}(\widetilde{U}). \tag{37}$$

From the above equation, together with (35), we conclude that

$$0 \in \partial_z f(z^*) + N_Z(z^*) + \partial_z Q(z^*, \xi^*), \tag{38}$$

i.e., $z^*$ is an optimal solution to the one-scenario problem $\min_{z \in Z} f(z) + Q(z, \xi^*)$. Clearly, if the latter problem has a unique optimal solution, then it must coincide with $z^*$. □

Note that the uniqueness of optimal solutions of $\min_{z \in Z} f(z) + Q(z, \xi^*)$ is assured for example when $f$ is *strictly* convex—a property that can be enforced by adding a regularization term if needed.

Theorem 1 shows the existence of an optimal scenario. As discussed earlier, it is clear that such a result is not of immediate use as it depends on optimal solutions of the same problem. In next sections, however, we will see how Theorem 1 plays a fundamental role in the development of approximating policies.

Still, some conclusions can be drawn from Theorem 1. For instance, the special case when $T$ is random with $\mathbb{E}[T] = 0$ leads to a somewhat surprising result, which to the best of our knowledge has not been observed in the literature. Essentially, it says that in that setting the second stage is irrelevant for the first-stage decision. While the result is obvious when $T$ is constant equal to zero—since in that case there is no link between the first and second stages—it is interesting to see that such a property remains valid when $T$ is random with $\mathbb{E}[T] = 0$.

**Corollary 1.** *Suppose Assumptions 1-2 hold. Suppose also that the random matrix $T$ is such that $\mathbb{E}[T] = 0$. Then, if the first-stage problem $\min_{z \in Z} f(z)$ has a unique optimal solution $z^*$, then $z^*$ is an optimal solution of* (2)-(5).

*Proof.* By Theorem 1, an optimal scenario $\xi^*$ is given by $T^* = \mathbb{E}[T] = 0$ and $h^* = T^* z^* = 0$. It follows that $Q(z, \xi^*)$ does not depend on $z$ and hence the optimal solutions of the one-scenario problem $\min_{z \in Z} f(z) + Q(z, \xi^*)$ are the same as those of the problem $\min_{z \in Z} f(z)$. Again, Theorem 1 ensures that if the latter problem has a unique optimal solution, then such a solution is optimal for the two-stage problem (2)-(5). □

Theorem 1 also provides an interesting support for a practice, observed in the context of energy system operators, of solving operational planning problems with only one scenario of demand instead of modeling the entire distribution. As discussed in Dias-Garcia et al. (2024), there are reasons for such a choice. System operators then compute pointwise forecasts of the demand and add a heuristic bias to it. Theorem 1 shows that, as long as the *optimal* bias is added, the practice is actually correct, in the sense that the final result is the same as though the full distribution were used. Indeed, the numerical results in Dias-Garcia et al. (2024)—using a method that actually motivated the development of the present paper—show that the use of a carefully computed bias in that setting yields very satisfactory results compared to practical benchmarks. Our results in Section 7 also corroborate that idea, now with the theoretical support provided by Theorem 1 as well as the results in Sections 4 and 5.

## 3.1 Finding an optimal scenario

We discuss now a possible way to search for an optimal scenario. Under the assumptions of Theorem 1, we know that there exists $\xi^*$ such that the optimal solution of the *one-scenario* problem $\min_{z \in Z} f(z) + Q(z, \xi^*)$ coincides with the optimal solution $z^*$ of (2)-(5). Thus, we can search directly for $\xi^*$ by solving the bilevel problem

$$\min_{\hat{\xi} \in \mathbb{R}^m} \mathbb{E}\left[G(z_D^*(\hat{\xi}), \xi)\right] \tag{39}$$

$$\text{s.t.} \quad z_D^*(\hat{\xi}) = \underset{z \in Z}{\operatorname{argmin}} \ G(z, \hat{\xi}) \tag{40}$$

where $G(z, \xi) = f(z) + Q(z, \xi)$. It is clear that any optimal scenario $\xi^*$ for (2)-(5) solves (39)-(40), since $z_D^*(\xi^*) = z^*$. Solving (39)-(40), of course, is not an easy task, since typically it is a non-convex problem. However, it is interesting to note that the distribution of $\xi$ is only used to *evaluate* a candidate solution $z_D^*(\hat{\xi})$. A consequence of this fact is that, when sampling is used to approximate the expected value in (39), the complexity of (39)-(40) grows *linearly* with the sample size. This property can be contrasted with solving (2)-(5) by using a Sample Average Approximation (SAA) approach and solving the resulting sampled problem as a linear program—the complexity of which grows exponentially with the sample size, as we add new variables and constraints for each sample. We shall see a numerical illustration of such a phenomenon in Section 7.

## 4 Optimal solution mappings

As discussed in Section 1, we apply the one-scenario result in Theorem 1 to the case of stochastic optimization problems with contextual information formulated in (9). Our goal is to derive approximations to the function that maps a contextual information $x$ to an optimal solution of (9). We start by formalizing the notion of optimal solution mappings discussed in the introduction. Let $z_S^* : \mathbb{R}^s \mapsto \mathbb{R}^n$ be a mapping defined such that $z_S^*(x)$ is an optimal solution of (9) for $x \in \mathbb{R}^s$ (we assume that (9) does indeed have an optimal solution for any $x \in \mathbb{R}^s$). Such a definition conveys the fact that, in practice, the decision maker will likely not be interested in solving the problem for a particular value of $x$; rather, the decision maker would like to have a *policy* that yields (or, more realistically, approximates) the optimal solution of (9) for any given $x$. In that sense, a policy $\pi$ is a mapping $\pi : \mathbb{R}^s \mapsto \mathbb{R}^n$ such that $\pi(x)$ is feasible (i.e., $\pi(x) \in Z$) for any $x \in \mathbb{R}^s$). Such a policy, of course, must be constructed from available data.

A natural question that arises then is, how to evaluate one such policy? One way to do this is by computing the *out-of-sample* performance of that policy. That is, given a dataset $(x_1, \xi_1), \ldots, (x_t, \xi_t)$, we use part of the dataset as training data to construct the policy $\pi$, and then use the remaining dataset (the testing data) to evaluate out-of-sample performance. This is a standard procedure but it illustrates the need for policies rather than seeking for just an optimal solution for a given $x$. In fact, if the contextual information of interest $x$ does not appear in the testing dataset, or if there are only a handful of observations with that $x$, we cannot really evaluate the objective function of (9). Such an issue of course is also present in the training data, but in that case building approximations of the conditional distribution of $\xi$ with respect to $x$ constitutes one way to generate a policy, as seen in Section 1. In the out-of-sample evaluation we cannot resort to such approximations, otherwise we would be distorting the actual value of the objective function of (9).

The above discussion suggests that what we are really interested in is the average performance of a policy $\pi$ over the set of features, i.e.,

$$\mathcal{P}(\pi) := \int_{\mathcal{X}} \mathbb{E}_\xi \left[ G(\pi(x), \xi) \,|\, X = x \right] F_X(dx), \tag{41}$$

where $F_X$ is the distribution of the features $X$. Note that we can rewrite (41) as

$$\begin{aligned} \mathcal{P}(\pi) &= \mathbb{E}_X \left[ \mathbb{E}_\xi \left[ G(\pi(X), \xi) \,|\, X \right] \right] \\ &= \mathbb{E}_{X,\xi} \left[ G(\pi(X), \xi) \right]. \end{aligned} \tag{42}$$

The search for the best policy can then be formulated as

$$\min_{\pi \in \Pi} \mathcal{P}(\pi), \tag{43}$$

where $\Pi$ is the set of mappings $\mathcal{X} \mapsto Z$. The notion of a policy in the above context is discussed in Sadana et al. (2024).

An example of a mapping $\pi \in \Pi$ is given by optimal solutions of (1), i.e., the solution of the problem that ignores the contextual information:

$$\pi_{NF}(x) \in \underset{z \in Z}{\operatorname{argmin}} \ \mathbb{E}_\xi \left[ G(z, \xi) \right] \quad \forall x \in \mathcal{X}, \tag{44}$$

where the subscript $NF$ stands for "no features". As $\pi_{NF}(\cdot) \in \Pi$, we have that

$$\mathbb{E}_{X,\xi} \left[ G(\pi_{NF}(X), \xi) \right] \geq \min_{\pi \in \Pi} \mathbb{E}_{X,\xi} \left[ G(\pi(X), \xi) \right],$$

which conveys the intuitive notion that ignoring the contextual information may lead to sub-optimal solutions.

Another example, of course, is given by optimal solutions of (9):

$$\pi_S(x) \in \underset{z \in Z}{\operatorname{argmin}} \ \mathbb{E}_\xi \left[ G(z, \xi) \,|\, X = x \right] \quad \forall x \in \mathcal{X}. \tag{45}$$

Proposition 1 below shows that $\pi_S$ in fact characterizes the optimal solutions of (43). Although the result is mentioned in Sadana et al. (2024) as a consequence of an interchangeability property from Rockafellar and Wets (1998), we present a proof here to make the paper self-contained and to add some intermediate steps to that argument.

**Proposition 1.** *Suppose Assumptions 1-2 hold. Then, the mapping $\pi_S$ defined in (45) solves (43). Moreover, if $\pi^*$ solves (43), then $\pi^*(x) \in \underset{z \in Z}{\operatorname{argmin}} \ \mathbb{E}_\xi \left[ G(z, \xi) \,|\, X = x \right]$ for all $x \in \mathcal{X}$ except perhaps on a set of $F_X$-measure zero.*

*Proof.* Consider the function $f : \mathbb{R}^n \times \mathbb{R}^s \mapsto \mathbb{R}$ defined as

$$f(z, x) := \mathbb{E}_\xi \left[ G(z, \xi) \,|\, X = x \right].$$

Then, from (45) we have that $\pi_S(x) \in \underset{z \in Z}{\operatorname{argmin}} \ f(z, x)$. Moreover, we have that

$$\begin{aligned}
\min_{\pi \in \Pi} \mathcal{P}(\pi) &= \min_{\pi \in \Pi} \mathbb{E}_{X,\xi} \left[ G(\pi(X), \xi) \right] \\
&= \min_{\pi \in \Pi} \int_\mathcal{X} \mathbb{E}_\xi \left[ G(\pi(x), \xi) \,|\, X = x \right] F_X(dx) \\
&= \min_{\pi \in \Pi} \int_\mathcal{X} f(\pi(x), x) \, F_X(dx) \tag{46} \\
&= \int_\mathcal{X} \min_{z \in Z} f(z, x) \, F_X(dx) \tag{47} \\
&= \int_\mathcal{X} f(\pi_S(x), x) \, F_X(dx) \\
&= \int_\mathcal{X} \left( \mathbb{E}_\xi \left[ G(\pi_S(x), \xi) \,|\, X = x \right] \right) F_X(dx) \\
&= \mathbb{E}_{X,\xi} \left[ G(\pi_S(X), \xi) \right] = \mathcal{P}(\pi_S),
\end{aligned}$$

where the equality in (47) follows from Theorem 14.60 in Rockafellar and Wets (1998). Thus, $\pi_S$ solves (43). The converse statement also follows from the same Theorem 14.60, which ensures—under a finiteness assumption that holds under Assumptions 1-2—that any $\pi$ that solves (46) must coincide with the mapping defined by an optimal solution of the inner problem in (47) for $F_X$-almost all $x \in \mathcal{X}$. □

The result in Proposition 1 is intuitive—the mapping that yields an optimal solution for each $x$ also yields an optimal solution on the average, and vice-versa. Approximating that mapping, however, is a difficult task. By using Theorem 1, however, we can obtain a stronger result in the case of two-stage problems of the form (2)-(5) under an additional assumption stated below.

**Assumption 3.** *The one-scenario problem $\min_{z \in Z} f(z) + Q(z, \xi)$ has a unique optimal solution for each value of $\xi \in \mathbb{R}^m$.*

As discussed earlier, the assumption holds for example when $f$ is strictly convex. When $f$ is linear, the assumption can be enforced by adding a regularization term, or by perturbing the coefficients $c$ and $q$ as discussed in Dias-Garcia et al. (2024).

**Proposition 2.** *Consider the mapping $\pi_S$ defined in (45), and suppose that $G(z, \xi) = f(z) + Q(z, \xi)$ with $Q$ defined in (3)-(5). Also, suppose Assumptions 1-3 hold. For each $x \in \mathcal{X}$, define $\xi_x^*$ as an optimal scenario to (9)—the existence of which is ensured by Theorem 1. Then, by defining $\pi_D(x)$ as the unique optimal solution of one-scenario problem $\min_{z \in Z} f(z) + Q(z, \xi_x^*)$ we have that the mapping $\pi_D$ is the unique solution to (43), except perhaps on a set of $F_X$-measure zero.*

*Proof.* The proof follows the same steps as in the proof of Proposition 1, noting that under the assumptions of the proposition we have that $\pi_D(x) = \pi_S(x)$ by virtue of Theorem 1. □

Proposition 2 has an important consequence: since $\pi_D(x)$ can be obtained simply by solving (6)-(8) with $\xi_x^*$ in place of $\xi^*$, to approximate the optimal mapping $\pi^*$ that solves (43) it suffices to approximate an optimal scenario $\xi_x^*$ for each $x$. This leads to the notion of "optimal pointwise forecasts". We discuss this topic in detail in the next section.

## 5 Optimal pointwise forecasts

Proposition 2 demonstrates that, when $G(z, \xi) = f(z) + Q(z, \xi)$ with $Q$ defined in (3)-(5)—which we will assume throughout this section—the search for policies that approximately solve (43) can be reduced to the search of policies that approximate the mapping $\pi_D(\cdot)$. One way to derive approximations to $\pi_D(x)$ is by approximating an optimal scenario $\xi_x^*$ with another function of $x$. To see this, notice that Assumptions 1-3 imply that the function $z_D^*(\cdot)$ defined as

$$z_D^*(\xi) := \operatorname*{argmin}_{z \in Z} f(z) + Q(z, \xi) \tag{48}$$

is *continuous*; see, e.g., Corollaries 8.1 and 9.1 in Hogan (1973). Since $\pi_D(x) = z_D^*(\xi_x^*)$, it follows that if $\psi(x)$ is a function such that $\psi(x) \approx \xi_x^*$, then $\widetilde{\pi}_D(x) := z_D^*(\psi(x))$ yields an approximation to $\pi_D(x)$.

In light of the above discussion, our goal is to define a data-driven approximation to $\xi_x^*$, constructed from observations $(x_1, \xi_1), \ldots, (x_N, \xi_N)$. To do so, we shall use a cost-based measure of forecast error, rather than a purely statistical measure. As discussed in Section 1, there is a growing body of literature on forecast models that are tailored to the optimization problem where such forecast is used, as it has been demonstrated that such models may indeed lead to better solutions. In our case, this forecast error is measured as the cost difference between using the forecast-driven decision and an ideal (i.e. *ex-post*) decision. This approach leads naturally to a bilevel optimization formulation. To see that, define $\Psi(\theta, x)$ as a *pointwise forecast function* of $\xi$ as a function of $x$, parameterized by $\theta$. Then, the optimal parameter $\theta_N^*$ solves the bilevel problem

$$\min_{\theta \in \Theta} \frac{1}{N} \sum_{n=1}^{N} \left| G(z_D^*(\hat{\xi}_n), \xi_n) - G(z_D^*(\xi_n), \xi_n) \right| \tag{49}$$

$$\text{s.t.} \quad \hat{\xi}_n = \Psi(\theta, x_n), \quad n = 1, \ldots, N \tag{50}$$

$$z_D^*(\hat{\xi}_n) = \operatorname*{argmin}_{z \in Z} \ G(z, \hat{\xi}_n), \quad n = 1, \ldots, N \tag{51}$$

$$z_D^*(\xi_n) = \operatorname*{argmin}_{z \in Z} \ G(z, \xi_n), \quad n = 1, \ldots, N, \tag{52}$$

and we let

$$\psi(x) := \Psi(\theta_N^*, x). \tag{53}$$

As seen from the above formulation, Model (49)-(52) measures the error between the forecast $\hat{\xi}_n$ and the observed data $\xi_n$ in terms of the cost of using the respective optimal solutions. More specifically, given the scenario realization $\xi_n$, $z_D^*(\xi_n)$ given by (52) would have been the best possible decision for that scenario, thereby realizing a cost of $G(z_D^*(\xi_n), \xi_n)$. On the other hand, since we only have a forecast $\hat{\xi}_n = \Psi(\theta, x_n)$ of $\xi_n$, we compute the decision $z_D^*(\hat{\xi}_n)$ given by (51), which after the realization of $\xi_n$ incurs a cost of $G(z_D^*(\hat{\xi}_n), \xi_n)$. Thus, we want to bias the forecast function $\Psi(\theta, \cdot)$ so that it minimizes the *ex-post* total forecast regret $\frac{1}{N} \sum_{n=1}^{N} |G(z_D^*(\hat{\xi}_n), \xi_n) - G(z_D^*(\xi_n), \xi_n)|$.

Note that, by definition of $z_D^*$, we must have $G(z_D^*(\hat{\xi}_n), \xi_n) \geq G(z_D^*(\xi_n), \xi_n)$ so we can remove the absolute value in (49). Moreover, it is clear that $G(z_D^*(\xi_n), \xi_n)$ is constant for the optimization problem in $\theta$. It follows that (49)-(52) can be equivalently written as follows, as in Dias-Garcia et al. (2024) and Muñoz et al. (2022):

$$\min_{\theta \in \Theta} \frac{1}{N} \sum_{n=1}^{N} G(z_D^*(\hat{\xi}_n), \xi_n) \tag{54}$$

$$\text{s.t.} \quad \hat{\xi}_n = \Psi(\theta, x_n), \quad n = 1, \ldots, N \tag{55}$$

$$z_D^*(\hat{\xi}_n) = \operatorname*{argmin}_{z \in Z} \ G(z, \hat{\xi}_n), \quad n = 1, \ldots, N. \tag{56}$$

It is interesting to notice that, in the case where there is no contextual information, problem (54)-(56) reduces to the bilevel model (39)-(40) introduced in Section 3.1, with the mapping $\Psi(\theta) := \theta$ and the expectation replaced by a sample average (note that in this case (55)-(56) are identical for all $n$ so we only to need to write them once). The importance of Model (54)-(56) lies in Theorems 2 and 3 below, which show that the model yields in the limit the *best possible policy* among those based on forecasts parameterized by $\theta$. Theorem 2 is shown in Dias-Garcia et al. (2024) for the case where $G(z, \xi) = c^\top z + Q(z, \xi)$, by using a similar argument to—but with weaker assumptions than—than that used by Shapiro and Xu (2008) for more general optimization problems with equilibrium constraints. A closer look at that proof shows that the main requirements are the continuity of the function $z_D^*(\cdot)$—which follows from Assumptions 1-3, as discussed above—and integrability of $G(z_D^*(\hat{\xi}_n), \xi_n)$ as a function of $\theta$, which is ensured by conditions (iv)-(vi) of the theorem. It follows that the result can be easily extended to the case in (6)-(8) where $f$ is a convex function and $Z$ is a convex set; we present the theorem here for completeness and to state the result in our notation.

**Theorem 2.** *Consider Model* (54)-(56). *Suppose that (i) Assumptions 1-3 hold, (ii) the forecasting function*

$\Psi(\cdot, \cdot)$ *is continuous in both arguments, (iii) the data process* $(X_1, \xi_1), \ldots, (X_N, \xi_N)$ *is independent and identically distributed (i.i.d.), (iv) the random variable $\xi$ is integrable, (v) the feasibility set $Z$ is bounded, and (vi) the set $\Theta$ is compact and non-empty. Then, with probability 1,*

$$\lim_{N \to \infty} d(\theta_N^*, S^*) = 0, \tag{57}$$

*where $d$ is the Euclidean distance from a point to a set and $S^*$ is defined as*

$$S^* = \operatorname*{argmin}_{\theta \in \Theta} \ \mathcal{P}(\widetilde{\pi}_D^\theta), \tag{58}$$

*with*

$$\widetilde{\pi}_D^\theta(x) \ := \ z_D^*(\Psi(\theta, x)), \tag{59}$$

$\mathcal{P}$ *defined in* (42), *and* $z_D^*$ *defined in* (48).

Note that the assumption that the data process $(X_1, \xi_1), \ldots, (X_N, \xi_N)$ is i.i.d. is a common assumption in the literature. Also, in Dias-Garcia et al. (2024) the result is extended to the case where $X_n$ is a (measurable) function of $\xi_1, \ldots, \xi_{n-1}$, and the data process generating $\{\xi_n\}_{n=1}^\infty$ is a *stationary ergodic* time series. Such a situation covers the case where the contextual information actually consists of previous observations, or some function thereof.

Theorem 2 assumes that the forecasting function $\Psi(\cdot, \cdot)$ is continuous in both arguments. Such an assumption covers many cases of interest, such as when the forecast value is an affine function of the contextual information (as in regression), or more generally when $\Psi$ is given by a neural network built upon continuous functions such as ReLu.

For some methods, however, continuity does not hold; this is the case for example of classification trees, or more generally of classification and regression trees (CART). When the continuity assumption does not hold, convergence can still be achieved as long as the parameter $\theta$ takes on only finitely many values. In that case, Theorem 3 below provides a reasonable alternative to Theorem 2.

**Theorem 3.** *Consider Model* (54)-(56). *Suppose that (i) Assumptions 1-3 hold, (ii) the data process $(X_1, \xi_1), \ldots, (X_N, \xi_N)$ is independent and identically distributed (i.i.d.), (iv) the random variable $\xi$ is integrable, (v) the feasibility set $Z$ is bounded, and (vi) the set $\Theta$ is* finite *and non-empty. Then, the event*

$$\theta_N^* \in S^* \tag{60}$$

*happens w.p.1 for $N$ large enough, where $S^*$ is defined in* (58). *If in addition, the support $\Xi$ of $\xi$ is bounded, then the convergence occurs* exponentially fast, *in the sense that there exist positive constants $K$ and $\beta$ such that*

$$P(\theta_N^* \notin S^*) \ \leq \ K e^{-N\beta}. \tag{61}$$

*Proof.* For the first claim, it suffices to show that the conditions of Proposition 2.1 in Kleywegt, Shapiro, and Homem-de-Mello (2002) hold in this case. As in that paper, Model (54)-(56) is a discrete stochastic optimization problem with finite feasibility set $\Theta$. Moreover, as in the proof of Theorem 2, under the assumptions of the theorem the random variable $G(z_D^*(\Psi(\theta, X_n), \xi_n))$ is integrable for all $\theta$. Proposition 2.1 in Kleywegt et al. (2002) then ensures that (60) holds w.p.1 for $N$ large enough.

The second claim follows from the fact that, under the assumption on boundedness of $\Xi$, the random variable $G(z, \xi)$ is bounded for all $z \in Z$ and hence Proposition 2.2 in Kleywegt et al. (2002) ensures that (61) holds. Note that the assumption on boundedness of $\Xi$ can be relaxed to finiteness of a certain moment generating function in a neighborhood of zero; we refer to Kleywegt et al. (2002) for details. $\square$

The finiteness assumption on $\Theta$ actually holds when standard CART is used as a forecast method, as long as the contextual information $X$ takes on only *finitely many* values. Indeed, in that case the forecast function $\Psi(\theta, x)$ yields the average of the observations $\xi_n$ in each leaf of the tree specified by $\theta$. Such specification consists of the branching order of the components of $x$, together with the threshold associated with each branching. When the contextual information $X$ takes on only finitely many values, it is clear that the set of possible thresholds can be reduced to the set of values taken by $X$. Thus, in that case the set of possible trees —and hence the set of possible values of $\theta$—is finite. The finiteness property can also be seen from the mixed integer programming formulations for CART in Bertsimas and Dunn (2017) and Verwer and Zhang (2019), where it is shown that the branching order is modeled with a finite number of binary variables that depends on the depth of the tree, which is fixed *a priori* as a parameter of the method; the latter paper also shows that the thresholds can be modeled with finitely many binary variables, the number of which depend on the maximum number of distinct values for each feature.

## 5.1 Approximating the optimal forecast

The results in Theorems 2 and 3 show that the policy $\widetilde{\pi}_D^{\theta_N^*}$ defined in (59) converges to the best possible policy obtained with the forecast function $\Psi$. It remains to study how far the policies $\widetilde{\pi}_D^\theta$ are from the policy that solves (43). Since $\widetilde{\pi}_D^\theta(x) = z_D^*(\Psi(\theta, x))$, we see from Proposition 2 that the answer to that question lies in how well $\Psi(\theta_N^*, x)$ approximates an optimal scenario $\xi_x^*$. This is summarized in Theorem 4 below, which is proved for the case where the first-stage problem in (2) is linear:

**Theorem 4.** *Consider the case where the first-stage problem in (2) is linear, i.e., $f(z) = c^\top z$ and $Z$ is polyhedral. Suppose that there exist $\delta \geq 0$ and $\widehat{\theta} \in \Theta$ such that $\|\Psi(\widehat{\theta}, x) - \xi_x^*\| \leq \delta$ for all $x \in \mathcal{X}$. Then, under the assumptions of either Theorem 2 or Theorem 3, there exists a constant $K \geq 0$—which depends only on the parameters that define the function $G$—such that the policy $\widetilde{\pi}_D^{\theta_N^*}$ obtained from Model (54)-(56) satisfies*

$$\lim_{N \to \infty} \mathcal{P}\big(\widetilde{\pi}_D^{\theta_N^*}\big) - \min_{\pi \in \Pi} \mathcal{P}\big(\pi\big) \ \leq \ K\delta. \tag{62}$$

*Proof.* Under conditions of the theorem, the function $z_D^*$ defined in (48) (which we write here as $z_D^*(h, T)$) is the optimal solution of the linear program

$$\min_{z \in Z} c^\top z + q^\top y$$
$$\text{s.t.} \quad Tz + Wy = h$$
$$y \geq 0.$$

It follows from (Robinson, 1973, Corollary 3.1) that $z_D^*(h, T)$ is a Lipschitz function of $(h, T)$ with constant, say, $M_1$. Thus, by the assumption in the theorem there exist $\delta \geq 0$ and $\widehat{\theta}$ such that $\|\Psi(\widehat{\theta}, x) - \xi_x^*\| \leq \delta$ for all $x \in \mathcal{X}$, which implies that

$$\left\|\widetilde{\pi}_D^{\widehat{\theta}}(x) - \pi_D(x)\right\| \ = \ \left\|z_D^*(\Psi(\widehat{\theta}, x)) - z_D^*(\xi_x^*)\right\| \ \leq \ M_1\delta. \tag{63}$$

Next, note that we have, for sufficiently large $N$,

$$\mathcal{P}(\widetilde{\pi}_D^{\theta_N^*}) - \min_{\pi \in \Pi} \mathcal{P}(\pi) = \mathcal{P}(\widetilde{\pi}_D^{\theta_N^*}) - \mathcal{P}(\pi_D)$$
$$\leq \mathcal{P}(\widetilde{\pi}_D^{\widehat{\theta}}) - \mathcal{P}(\pi_D). \quad (64)$$

The inequality in (64) follows from the fact that under the assumptions of either Theorem 2 or Theorem 3 the sequence $\{\theta_N^*\}$ approaches the optimal set $S^*$ and so for $N$ large enough we have that $\mathcal{P}(\widetilde{\pi}_D^{\theta_N^*}) \leq \mathcal{P}(\widetilde{\pi}_D^{\widehat{\theta}})$. Moreover, from (30) we have that the subdifferential set of $G(z,\xi)$ is bounded for all $z \in Z$ and $\xi \in \Xi$ which in turn implies that $G(\cdot,\xi)$ is uniformly Lipschitz, i.e., there exists a constant $M_2 \geq 0$ such that

$$|G(z_1,\xi) - G(z_2,\xi)| \leq M_2 \|z_1 - z_2\| \quad \text{for all } z_1, z_2 \in Z \text{ and all } \xi \in \Xi. \quad (65)$$

Inequalities (63) and (65), together with (64) and definition (42) of $\mathcal{P}$, then imply that

$$\mathcal{P}(\widetilde{\pi}_D^{\theta_N^*}) - \min_{\pi \in \Pi} \mathcal{P}(\pi) \leq K\delta,$$

where $K := M_1 M_2 \delta$. $\square$

Theorem 4 materializes the notion of "optimal pointwise forecasts": as long as an optimal scenario $\xi_x^*$ (viewed as a function of $x$) can be approximated uniformly by some function $\Psi(\theta,\cdot)$, a *pointwise* forecast constructed with the parameters $\theta_N^*$ resulting from Model (54)-(56) will suffice, in the sense that the policy $\widetilde{\pi}_D^{\theta_N^*}(x)$ defined as $z_D^*(\Psi(\theta_N^*,x))$ will yield an approximate solution to (43). As indicated by (62), the quality of the latter approximation depends on how well $\Psi(\theta,\cdot)$ approximates $\xi_{(\cdot)}^*$. Note that, while the form of $\Psi$ must be specified in advance, the actual values of $\theta$ that make $\Psi(\theta,\cdot)$ approximate $\xi_{(\cdot)}^*$ need not be known in advance—in fact, they result from applying Model (54)-(56) as stated in Theorem 4. Thus, the more functions the set $\{\Psi(\theta,\cdot) : \theta \in \Theta\}$ contains, the better the approximation of the policy that solves (43).

Naturally, in order to be able to approximate $\xi_{(\cdot)}^*$ we need this function to have some properties. One such property is described in Proposition 3 below.

**Proposition 3.** *Suppose that $\mathcal{X}$ is a continuous set and that the conditional distribution of $\xi|x$ [2] is close to the distribution of $\xi|x'$ when $x$ is close to $x'$. More precisely, suppose that given $\varepsilon > 0$, there exists $\delta > 0$ such that*

$$\|x - x'\| < \delta \implies d_W(\xi|x, \xi|x') < \varepsilon, \quad (66)$$

*where $d_W$ denotes the Wasserstein distance between two distributions. Then, if the mapping $\pi_S(\cdot)$ defined in (45) is single-valued for each $x$, then $\pi_S(\cdot)$ is continuous and hence $\xi_{(\cdot)}^*$ constructed in Theorem 1 is continuous.*

*Proof.* The results follows from classical stability results for optimal solutions of stochastic programs. For instance, Corollary 14 in Römisch (2003) shows that, in the case of the two-stage stochastic program (2)-(5), there exists $\overline{\delta} > 0$ such that, if $P$ and $Q$ are two distributions such that $d_W(P,Q) < \overline{\delta}$, then we have that

$$Z^*(Q) \subseteq Z^*(P) + \varphi(L\, d_W(P,Q))\, \mathbb{B}$$

where $Z^*(Q)$ and $Z^*(P)$ are the optimal solution sets of (2)-(5) when the distribution of $\xi$ is respectively $Q$ and $P$, $\varphi(\cdot)$ is a certain increasing function that vanishes at zero, $L$ is a positive constant, and $\mathbb{B}$ is the

---
[2] Here we abuse the notation—by "the conditional distribution of $\xi|x$" we mean the conditional probability measure defined as $P_x(A) := P(\xi \in A \mid X = x)$ for Borel sets $A$.

Euclidean unit ball. Thus, if $Z^*(P)$ and $Z^*(Q)$ are singletons, it follows that given $\eta > 0$ sufficiently small, there exists $0 < \varepsilon < \bar{\delta}$ such that $\|Z^*(P) - Z^*(Q)\| < \eta$ whenever $d_W(P,Q) < \varepsilon$. By putting the conditional distributions of $\xi|x$ and $\xi|x'$ in place of $P$ and $Q$, from (66) we see that the condition $d_W(P,Q) < \varepsilon$ holds whenever $\|x - x'\| < \delta$. We conclude that $\pi_S(\cdot)$ is continuous. Since $\xi_x^* = T^*\pi_S(x)$, we see that in that case $\xi_{(\cdot)}^*$ is continuous as well. □

A particular case where condition (66) in Proposition 3 is satisfied is when the uncertainty $\xi$ can be written as a linear model of $x$, i.e.,
$$\xi = A + Bx + \epsilon,$$
where $A$ is a constant vector, $B$ is a matrix of coefficients, and $\epsilon$ has a multi-variate Normal distribution with mean 0 and covariance matrix $\Sigma$. That is, for fixed $x$ and $x'$ we have that
$$\xi|x \sim \text{Normal}(A + Bx, \Sigma), \quad \xi|x' \sim \text{Normal}(A + Bx', \Sigma)$$
which are close when $x$ is close to $x'$.

## 5.2 Solving the bilevel problem

We discuss now some methods to solve the bilevel problem (54)-(56) which is constructed from observations $(x_1, \xi_1), \ldots, (x_N, \xi_N)$ of the feature $X$ and the random variable $\xi$. Throughout this section, we will assume that the matrix $T$ in (4) is not random[3], so the random variable $\xi$ corresponds only to the right-hand side term $h$. Recalling that $G(z, \xi)$ is defined as $c^\top z + Q(z, \xi)$ (with $Q$ defined in (3)-(5)), we can write (54)-(56) as

$$\min_{\theta, (z_1, y_1^u), \ldots, (z_N, y_N^u)} \frac{1}{N} \sum_{n=1}^{N} c^\top z_n + q^\top y_n^u \tag{67}$$
$$\text{s.t.} \quad \theta \in \Theta \tag{68}$$
$$z_n \in Z, \; y_n^u \geq 0, \quad n = 1, \ldots, N \tag{69}$$
$$Tz_n + Wy_n^u = \xi_n, \quad n = 1, \ldots, N \tag{70}$$
$$(z_n, y_n^\ell) = \operatorname*{argmin}_{z \in Z, \; y \geq 0} \left\{ c^\top z + q^\top y : Tz + Wy = \Psi(\theta, x_n) \right\}, \quad n = 1, \ldots, N. \tag{71}$$

Note that we have different second-stage variables for the upper and lower level problems, denoted respectively $y_n^u$ and $y_n^\ell$, $n = 1, \ldots, N$. Also, since the lower level problems are just linear programs, we can write (67)-(71) as a single level problem using KKT conditions, as customary in the bilevel literature. For instance, assuming

---
[3]This assumption is included only to ease the notation, and to avoid generating further bilinear terms in the formulations discussed in this section. The assumption can be relaxed, at the expense of requiring solving a problem that is, in principle, harder.

for simplicity that the set $Z$ is the positive orthant, we have that (67)-(71) is equivalent to

$$\min_{\theta, \{(z_n, y_n^u, y_n^\ell, \lambda_n)\}_{n=1}^N} \frac{1}{N} \sum_{n=1}^N c^\top z_n + q^\top y_n^u \tag{72}$$

$$\text{s.t.} \quad \theta \in \Theta \tag{73}$$

$$z_n \geq 0, \quad n = 1, \ldots, N \tag{74}$$

$$T z_n + W y_n^u = \xi_n, \ y_n^u \geq 0, \quad n = 1, \ldots, N \tag{75}$$

$$T z_n + W y_n^\ell = \Psi(\theta, x_n), \ y_n^\ell \geq 0, \quad n = 1, \ldots, N \tag{76}$$

$$T^\top \lambda_n \leq c, \quad n = 1, \ldots, N \tag{77}$$

$$W^\top \lambda_n \leq q, \quad n = 1, \ldots, N \tag{78}$$

$$c^\top z_n + q^\top y_n^\ell - \Psi(\theta, x_n)^\top \lambda_n \leq 0, \quad n = 1, \ldots, N. \tag{79}$$

In the above formulation, (73)-(75) are the upper level constraints, (76) represents the lower level primal constraints, (77)-(78) are the lower level dual constraints (with corresponding multipliers $\lambda_n$), and (79) imposes strong duality on the lower level, which is equivalent to writing the complementarity constraints.

Problem (72)-(79) is in general hard to solve; constraint qualifications typically required by nonlinear optimization algorithms do not hold, and the model becomes harder especially when the forecast function $\Psi(\theta, x)$ is a complicated function of $\theta$. A situation where the above model can be solved reasonably efficiently is when (i) the set $\Theta$ is polyhedral, and (ii) the function $\Psi(\theta, x)$ is linear in $\theta$, i.e. $\Psi(\theta, x) = B_x \theta$ for some matrix $B_x$. The latter condition corresponds to an "application-driven regression" whereby, instead of measuring the error with quadratic loss functions as in standard regression, we measure it using the objective function of the problem. In that case, the problem given by (72)-(79) is a linear program except for (79), which contains the bilinear term $\theta^\top B_{x_n}^\top \lambda_n$. Some alternatives to address that issue include using binary variables to eliminate the bilinear term (or to model the equivalent complementarity constraints), or relaxing (79) by putting it in the objective function. The latter can be solved with a penalty alternating direction method; we refer to Kleinert and Schmidt (2021) for further discussion and comparisons between these approaches.

When $\Theta$ is finite, as in Theorem 3, it may be possible to solve (72)-(79) as a mixed integer program, depending again on the form of $\Psi(\theta, x)$ as a function of $\theta$. A similar approach is discussed in Muñoz et al. (2022), where big-M constraints are used for complementarity slackness conditions. That work also present regularization techniques, in which these complementary conditions are relaxed. The authors argue that interior-point methods fail to obtain even a local optimal solution to the bilevel problem due to complementarity constraints. To overcome this difficulty, they present a regularization approach, in which the slackness condition constraint may be violated by at most a small tolerance parameter $\varepsilon > 0$. As $\varepsilon \to 0$, the approach aims to converge to a local optimum. This regularization allows the bilevel problem to be reformulated and solved through the KKT conditions. Designing efficient algorithms for solving such bilevel problems remains an open direction for future research.

## 5.3 The case with functional dependencies

We end this section by noticing that the optimal scenario constructed in Theorem 1 is defined according to the constraints of the second-stage problem rather than by the random variables present in those constraints. This is an important distinction, as it implies that dependencies between random variables in different constraints are not necessarily respected by the optimal scenario.

We illustrate this issue with a simple example. Consider the standard newsvendor model discussed in Section 2, which was formulated as a two-stage model in (13)-(16). Suppose we add the constraint $z \geq 2D$ to the second stage, which can be modeled as $y^d = -2D + z$, $y^d \geq 0$. In the notation of (3)-(5), we have $y = [y^+, y^-, y^d]^\top$, $W = \begin{bmatrix} 1 & -1 & 0 \\ 0 & 0 & 1 \end{bmatrix}$, $h = [-D, -2D]^\top$, and $T = [-1, -1]^\top$ (to ensure that the complete recourse assumption holds, we can add a slack variable $y^s \geq 0$ with high cost, such that $y^d - y^s = -2D + z$). Again, let $z^*$ be the optimal solution of the stochastic model. Then, the optimal scenario given by Theorem 1 is $h^* := Tz^* = [-z^*, -z^*]^\top$, so we see that the functional dependence between the components of $h$ (i.e, that the second one is twice the first one) is not respected. However, it is easy to check that solving the one-scenario problem with $h^*$ in place of $h$ will indeed yield the same solution $z^*$ as the stochastic model. Thus, as discussed in Section 5.1, as long as one is able to approximate $h^*$ with the forecast function $\Psi(\cdot, \cdot)$, an approximation of the optimal solution $z^*$ can still be obtained. Later, at the end of Section 7.3—and also in Appendix E—we will comment on the effect of this restriction on one of our case studies.

## 6 Solution Methods

In this section, we present different approaches to find policies $\pi \in \Pi$ to solve the problem (9). We begin by describing application-driven (AD) methods, which consider the structure of the problem and find a function $\Psi(\theta, x)$ that defines the policy $\tilde{\pi}_D^\theta(x) = z_D^*(\Psi(\theta, x))$. Then, we present benchmark methods, where we include standard predict-then-optimize methods and conditional sampling methods. We summarize all the methods in Table 1.

### 6.1 Application-Driven Forecasts Methods

For AD methods, we consider two $\Psi(\theta, x)$ forecast functions, linear regression and a method based on regression trees known as M5. In both cases, we seek to determine the optimal parameters for each function by solving the bilevel problem (54)-(56). We also propose a heuristic solution method.

As discussed earlier, problem (54)-(56) seeks a parameterization of the function $\Psi(\theta, x)$ such that the $z$ decisions obtained from this prediction function minimize the average cost obtained by considering observations $\xi_n$. Let $\theta_N^*$ be the parameterization found by solving the bilevel problem. The policy is given by,

$$\tilde{\pi}_D^{\theta_N^*}(x) \in \underset{z \in Z}{\operatorname{argmin}}\ G(z, \Psi(\theta_N^*, x)). \tag{80}$$

In general, the bilevel problem (54)-(56) is difficult to solve, including the problems presented in our computational study. For this reason, we consider heuristic methods to solve it, in particular, the implementation of the Meta algorithm and the Nelder-Mead method presented in Dias-Garcia et al. (2024). In this work, this heuristic method yields high-quality solutions in short computational running times. The Meta algorithm is

| Type | Method | Description |
|---|---|---|
| Application-driven methods | AD — Application-driven version of linear regression | For feature $x \in \mathcal{X}$, $\tilde{\pi}_N^{AD}(x)$ is the solution in $Z$ of $\min_{z \in Z} G\left(z, \Psi^{LS}(\theta_N^*, x)\right)$, with $\theta_N^*$ a solution of the bilevel problem (54)-(56) trained with observations $(x_1, \xi_1)$, …, $(x_N, \xi_N)$, using an affine function $\Psi(\theta, x) = \theta_0 + \theta^\top x$. |
| | M5+AD — Application-driven version of M5 forecasting | For feature $x \in \mathcal{X}$, $\tilde{\pi}_N^{M5+AD}(x)$ is the solution in $Z$ of $\min_{z \in Z} G\left(z, \Psi^{M5}(\theta_N^*, x)\right)$, with $\theta_N^*$ a solution of the bilevel problem (54)-(56) trained with observations $(x_1, \xi_1)$, …, $(x_N, \xi_N)$, using an M5 regressor $\Psi^{M5}$ as $\Psi$, i.e. a different affine function $\Psi$ is defined on each leaf of the tree. |
| Predict-then-optimize methods | LS — Least squares regression | For feature $x \in \mathcal{X}$, $\tilde{\pi}_N^{LS}(x)$ is the solution in $Z$ of $\min_{z \in Z} G\left(z, \Psi^{LS}(\hat{\theta}, x)\right)$, with $\Psi^{LS}(\hat{\theta}, \cdot)$ the least squares regressor trained with observations $(x_1, \xi_1)$, …, $(x_N, \xi_N)$. |
| | CART — classification and regression tree | For feature $x \in \mathcal{X}$, $\tilde{\pi}_N^{CART}(x)$ is the solution in $Z$ of $\min_{z \in Z} G\left(z, \Psi^{CART}(\hat{\theta}, x)\right)$, with $\Psi^{CART}(\hat{\theta}, \cdot)$ a CART regressor trained with observations $(x_1, \xi_1)$, …, $(x_N, \xi_N)$. |
| Conditional distribution methods | SAA — Sample Average Approximation (SAA) with no features | Regardless of the feature $x \in \mathcal{X}$, $\tilde{\pi}_N^{SAA}$ is the solution in $Z$ of the SAA problem $\min_{z \in Z} \frac{1}{N} \sum_{n=1}^{N} G(z, \xi_n)$. |
| | KNN — SAA with k-nearest neighbours | For feature $x \in \mathcal{X}$, $\tilde{\pi}_N^{KNN}(x)$ is the solution in $Z$ of $\min_{z \in Z} \frac{1}{|\mathcal{N}(x)|} \sum_{n \in \mathcal{N}(x)} G(z, \xi_n)$, with $\mathcal{N}(x)$ the indices of the $k$ data points $x_1$, …, $x_N$ nearest to $x$. |
| | ER-SAA — empirical residuals-based SAA | For feature $x \in \mathcal{X}$, $\tilde{\pi}_N^{ER-SAA}(x)$ is the solution in $Z$ of $\min_{z \in Z} \frac{1}{N} \sum_{n=1}^{N} G\left(z, \Psi^{LS}(\hat{\theta}, x) + \varepsilon_n\right)$, where $\varepsilon_n$ are forecast errors $\varepsilon_n := \xi_n - \Psi^{LS}(\hat{\theta}, x_n)$, with $\Psi^{LS}(\hat{\theta}, \cdot)$ the least squares regressor trained with observations $(x_1, \xi_1)$, …, $(x_N, \xi_N)$. |

Table 1: Summary of methods used in the computational experiments. Our proposed methods are the Application-driven methods, i.e., AD and M5+AD; the rest are extant methods, see Kannan et al. (2022).

as follows:

```
Output: Optimized θ
1  Initialize θ;
2  while Not converged do
3      Update θ;
4      foreach n ∈ {1, ..., N} do
5          Forecast: ξ̂_n ← Ψ(θ, x_n);
6          Plan Policy: z*_n ← argmin_{z ∈ Z} G(z, ξ̂_n);
7          Cost Assessment: cost_n ← G(z*_n, ξ_n);
8      end
9      Compute cost: cost(θ) ← ∑_{n=1}^{N} cost_n;
10 end
```

**Algorithm 1:** Meta algorithm.

This algorithm starts by initializing the vector $\theta \in \Theta$, which can be simply initialized with zero values. As long as the algorithm does not converge, the variables are updated. The convergence criterion is that the objective function decreases less than $\epsilon = 10^{-7}$ between two consecutive iterations. The update for $\theta$ aims to minimize the objective function, for which we need an optimization method. As in Dias-Garcia et al. (2024), we propose to use a derivative-free method, the Nelder-Mead approach. For each $n = 1, \ldots, N$, the forecast for observation $(x_n, \xi_n)$ is found using the current $\theta$ vector and the covariate vector $x_n$. Then, for the prediction $\hat{\xi}_n$, a plan $z_n$ is obtained. The cost of this plan is calculated using the observation $n$ value, $\xi_n$. Finally, the total cost for the current $\theta$, $cost(\theta)$, is calculated.

### 6.1.1 Linear Regression Prediction

For a linear regression, the parameterization of the $\Psi(\theta, x)$ function is given by,

$$\Psi(\theta, x) = \theta_0 + \sum_{i=1}^{s} \theta_i x_i. \tag{81}$$

The resulting policy from the linear function under the Application-Driven approach is denoted as $\tilde{\pi}_N^{AD}(x)$.

### 6.1.2 M5 Prediction

The M5 model is presented in Quinlan (1992). It is a generalization to classification and regression tree models (CART) that seeks to improve the prediction made in each region, using linear regression instead of an average value.

The CART method considers a partition of the feature space into $R$ regions, where for each partition $\ell \in \{1, \ldots, R\}$ a constant value is determined as a prediction (Hastie et al., 2009; Murphy, 2012). The partition is carried out by branching a decision tree and each region $\ell \in \{1, \ldots, R\}$ is a leaf of the tree. Tree structures allow to fit nonlinear functions, obtaining high-quality predictions. In general, the method for defining these partitions is heuristic; when a node is branched, the splitting variable $x_j$, $j = 1, \ldots, s$, and the cutoff value $r$ that has the greatest impact on reducing the variance of the training data are identified. This splitting variable and point are defined so that two half-planes are obtained:

$$\mathcal{R}_1(j,r) = \{n = 1, \ldots, N \mid x_n \leq r\} \text{ and } \mathcal{R}_2(j,r) = \{n = 1, \ldots, N \mid x_n > r\}. \tag{82}$$

We then seek the splitting variable $j$ and the split point $r$ that minimize

$$\min_{j,r} \left\{ \min_{\psi_1} \sum_{n \in \mathcal{R}_1(j,r)} (\xi_n - \psi_1)^2 + \min_{\psi_2} \sum_{n \in \mathcal{R}_2(j,r)} (\xi_n - \psi_2)^2 \right\}. \tag{83}$$

For given $j$ and $r$, the inner minimization is solved by

$$\psi_1 = \frac{1}{|\mathcal{R}_1(j,r)|} \sum_{n \in \mathcal{R}_1(j,r)} \xi_n \text{ and } \psi_2 = \frac{1}{|\mathcal{R}_2(j,s)|} \sum_{n \in \mathcal{R}_2(j,s)} \xi_n. \tag{84}$$

In practice, the splitting variable and the split point are determined by reviewing all the given observations and then determining the best pair $(j, r)$. This guarantees that the regions $\mathcal{R}_1(j, r)$ and $\mathcal{R}_2(j, r)$ are always non-empty.

The branching continues until the tree reaches a given maximum height or the number of data points in each leaf is less than or equal to a given value, or both. The prediction for a new vector of covariates $X$, corresponding to a region $\ell \in \{1, \ldots, R\}$, is the average value of the observations in that region, $\bar{\xi}_\ell$.

The M5 model considers an additional step for prediction given a vector $X$: using the observations in a region, a linear regression is estimated. We propose to perform this linear regression considering an AD approach. That is, once a tree structure is defined, where each region has a set of observations, for this set, we solve the bilevel problem (54)-(56) and consider a linear regression function for $\Psi(\theta, x)$, as in equation (81). For each leaf $\ell \in R$, we obtain a vector $\theta_\ell$, which takes into account the structure of the problem. The training of the tree follows the same heuristic used by CART. The resulting policy from using M5 under the Application-Driven approach is denoted as $\tilde{\pi}_N^{M5+AD}(x)$.

If the decision tree has a maximum tree size, with a fixed partition $\{1, \ldots, R\}$, constructed using i.i.d. observations, then the convergence of each parameter $\theta_{\ell N}$, $\ell \in \{1, \ldots, R\}$, is guaranteed, as shown in Proposition 4 below.

**Proposition 4.** *Consider a decision tree with a fixed partition $\{1, \ldots, R\}$, constructed using independent and identically distributed (i.i.d.) $N' < \infty$ observations. Also, consider Model (54)-(56). Suppose that (i) Assumptions 1-3 hold, (ii) the data process $(X_1, \xi_1), \ldots, (X_N, \xi_N)$ is i.i.d., (iii) the random variable $\xi$ is integrable, and (iv) the primal and dual feasibility sets $Z$ and $U$ are bounded. Then, with probability 1,*

$$\lim_{N \to \infty} d(\theta_{\ell N}^*, S_\ell^*) = 0, \tag{85}$$

*where $d$ is the Euclidean distance from a point to a set and $S_\ell^* = \underset{\theta \in \Theta}{\operatorname{argmin}} \ \mathbb{E}_{X, \xi} \left[ G\big(\tilde{\pi}_D^\theta(X), \xi\big) \mid (X, \xi) \in (\mathcal{X}_\ell, \Xi_\ell) \right]$, and $(\mathcal{X}_\ell, \Xi_\ell)$ is the support for region $\ell \in \{1, \ldots, R\}$.*

*Proof.* Let $\ell \in \{1, \ldots, R\}$ be a region of the decision tree. Let $h_\ell(N)$ be the number of elements in region $\ell$, with $N$ training data points. We want to show that $h_\ell(N)$ goes to infinity when $N$ goes to infinity. When the data are i.i.d., by the law of large numbers, we have that $\frac{h_\ell(N)}{N} \to p_\ell$ w.p. 1, where $p_\ell$ is the probability that an observation belongs to that region. We then have two cases:

1. $p_\ell > 0$: in this case, we have that for $N$ large, $h_\ell(N) \sim N p_\ell$ and hence $h_\ell(N) \to \infty$ w.p. 1.

2. $p_\ell = 0$: in this case, the probability that there are observations in the region $\ell$ is zero, *i.e.* $h_\ell(N) = 0$ for all $N$ w.p. 1. Since the region $\ell$ is not empty (at least one observation from $N'$ training data is in region $\ell$), this cannot happen.

We conclude that the number of elements in the region $\ell$ goes to infinity. Then, we apply Theorem 2, which guarantees convergence of the parameter $\theta_{\ell N}^*$. □

It is worth noting that data partitioning could also be carried out using alternative methods, such as k-nearest neighbors (KNN) (Morales et al., 2023). However, the resulting clusters would differ from those generated by the M5 model. Specifically, the M5 model partitions the set of observations with the goal of producing accurate predictions, whereas KNN relies solely on data similarity criteria to form groups, without directly accounting for predictive performance.

## 6.2 Benchmark Methods

In this section, we present benchmark methods, which we will use as a solution reference for our AD methods in our computational study. We start by describing predict-then-optimize methods (PO), commonly used in practice, and then conditional distribution methods (CD).

### 6.2.1 Predict then Optimize Framework

Predict-then-optimize (PO) methods consider two stages to find a policy. In the first stage, the $\theta \in \Theta$ parameters are found for the $\Psi(\theta, x)$ function, $x \in \mathcal{X}$, so that an error metric takes its smallest value.

This metric aims to measure the difference between the observed value and the prediction of the function $\Psi(\theta, x)$. In general, for observations $(x_1, \xi_1), \ldots, (x_N, \xi_N)$, the quadratic error, $\sum_{n=1}^{N} \|\xi_n - \Psi(\theta, x_n)\|^2$, is used as the error metric. Once the parameter vector $\hat{\theta}$ is determined, the policy $\pi \in \Pi$ is defined according to Equation (80).

For these PO methods we consider a linear function (Equation (81)) and a decision tree method. The tree method we use is CART, which, unlike M5, partitions the feature space and then makes predictions as averages of the observations. CART manages to achieve a better approximation of the data function, but it is not an AD method. In fact, to define the CART parameters (structure of the tree), the objective is to minimize the quadratic error. Similarly to M5, CART is usually constructed by means of a greedy heuristic (Hastie et al., 2009).

The resulting policy from the linear function is denoted as $\tilde{\pi}_N^{LS}(x)$ and the resulting policy from CART, $\tilde{\pi}_N^{CART}(x)$, $x \in \mathcal{X}$.

### 6.2.2 Conditional Distribution Methods

Conditional distribution methods (CD) seek to approximate the conditional expectation $\mathbb{E}_\xi[G(z, \xi)|X = x]$, $x \in \mathcal{X}$, in Equation (45), using a set of scenarios for each $x$. These methods also require two stages; in a first stage the scenarios for each $x$ must be determined, a set $\mathcal{N}(x)$, and then, in the second stage, the following problem is solved to find the conditional distribution policy $\tilde{\pi}(x)$:

$$\tilde{\pi}(x) \in \underset{z \in Z}{\mathrm{argmin}} \; \frac{1}{|\mathcal{N}(x)|} \sum_{n \in \mathcal{N}(x)} G(z, \xi_n), \qquad \forall x \in \mathcal{X}. \tag{86}$$

Note that here we are assuming that each scenario has the same weight (probability) in the objective function.

The simplest of these conditional methods ignores the context $x \in \mathcal{X}$, and the problem is simply solved by considering the observations $(x_1, \xi_1), \ldots, (x_N, \xi_N)$ given. For this, the policy is determined by solving these scenarios as done using the sample average approximation approach (SAA). We denote this approach as the SAA policy. Note that this method can also be seen as an approximation of the $\pi^{NF}$ policy of Equation (44). Given a set of observations, the policy obtained $\tilde{\pi}_N^{SAA}(x)$, $x \in \mathcal{X}$, is given by:

$$\tilde{\pi}_N^{SAA}(x) \in \underset{z \in Z}{\mathrm{argmin}} \; \frac{1}{N} \sum_{n=1}^{N} G(z, \xi_n), \qquad \forall x \in \mathcal{X}. \tag{87}$$

Other approaches to generate policies based on conditional scenarios are based on machine learning techniques. Using contextual information $x \in \mathcal{X}$ and predictive methods, scenarios for $x$ are generated. Then, a solution is found by solving an approximation problem as in (86). In Bertsimas and Kallus (2020), several prediction methods are proposed to determine these scenarios, such as $k$-nearest-neighbors (KNN), local linear regression (LOESS), CART and random forests (RF). As a benchmark for our computational experiments, we consider the KNN method, since in other computational studies it has also been a benchmark approach, showing good performance (Kannan et al., 2022). The KNN method determines the $k$ nearest neighbors for an $x \in \mathcal{X}$, $\mathcal{N}(x) = \left\{n = 1, \ldots, N \mid \sum_{j=1}^{N} \mathbb{I}[\|x - x_i\| \geq \|x - x_j\|] \leq k\right\}$, where $\|\cdot\|$ is a distance metric (*e.g.*, Euclidean distance). The resulting policy $\tilde{\pi}_N^{KNN}(x)$, is given by:

$$\tilde{\pi}_N^{KNN}(x) \in \underset{z \in Z}{\mathrm{argmin}} \; \frac{1}{|\mathcal{N}(x)|} \sum_{n \in \mathcal{N}(x)} G(z, \xi_n). \tag{88}$$

The last conditional distribution method we consider is the *empirical residuals-based SAA* (ER-SAA) one proposed in Kannan et al. (2022). In this method, we start by estimating $\theta$ by minimizing the squared error of the $\Psi(\theta, x)$ function of $\theta$ over the observations $(x_1, \xi_1), \ldots, (x_N, \xi_N)$, thus obtaining an optimal $\hat{\theta}$. This allows to obtain a point prediction $\Psi(\hat{\theta}, x)$ for each $x \in \mathcal{X}$. Next, for each data observation $n = 1, \ldots, N$, a forecast error (empirical residual) $\varepsilon_n := \xi_n - \Psi(\hat{\theta}, x_n)$ is computed. Thus, the conditional distribution (45) can be approximated, obtaining a $\tilde{\pi}_N^{ER-SAA}(x)$ policy given by:

$$\tilde{\pi}_N^{ER-SAA}(x) \in \operatorname*{argmin}_{z \in Z} \frac{1}{N} \sum_{n=1}^{N} G\left(z, \Psi(\hat{\theta}, x) + \varepsilon_n\right), \qquad x \in \mathcal{X}. \tag{89}$$

# 7 Computational Experiments

In this section, we present numerical experiments to support our theoretical results and to analyze the performance of the proposed methods. The goal of these experiments is to provide a case study illustrating our main result: solving stochastic optimization problems with fixed recourse and fixed costs can be achieved by solving a single-scenario problem. This scenario is generated using the bilevel model (39)-(40) in the case with no contextual information, and as an optimal pointwise forecast in the sense of Theorems 2 and 3 in the case with contextual information. In both cases, the corresponding bilevel models were solved using derivative-free methods—in the non-contextual case we used the library `BlackBoxOptim.jl`, whereas in the contextual case we used `BiLevelJUMP.jl`.

We begin by illustrating the one-scenario result in the case where there is no contextual information, using the newsvendor models described in Section 2. We then consider the contextual case by using problems from the literature with synthetic data, followed by our study using real data. The synthetic data problems, namely, the resource allocation problem and the shipment planning problem, are studied in Kannan et al. (2022), and our generation of data and parameters for these instances replicates the approach used in that work. In Section 7.2.1, we describe this instance generation method and model evaluation in detail. Recall that the methods used to solve the problems are summarized in Table 1.

## 7.1 The No-Contextual Case: Explicit Optimal Scenarios

We conduct two experiments to illustrate the use of the bilevel model (39)-(40) to find an optimal scenario. In the first one, we consider the newsvendor model with unreliable supplier discussed in Section 2, with parameters $c = \eta = 300$ and $p = \pi = 4000$, yielding $\phi = 0.928$. In order to compare with the analytical solution in (18), we used the same distributions that led to that formula, i.e., demand has uniform distribution in $(0, b)$ with $b = 100$ and the reliability factor $U$ has uniform distribution in $(0, 1)$. Since those distributions are continuous, we solved an SAA version of (39)-(40) with 1000 samples to determine an optimal scenario $(D^*, U^*)$. We then compared the corresponding solution given by $D^*/U^*$ with the solution of the SAA version of (18)-(23), computed with the same samples. Table 3 displays the results for five replications. As it can be seen from the table, in all instances model (39)-(40) produces (via optimal scenario) exactly the same solution as the SAA version of (18)-(23)—and both are reasonable approximations of the exact solution given by (18). Also, we see that in four of the five replications the optimal scenario found by (39)-(40) yields a demand scenario $D^*$ that lies outside of the support of the demand (which is $(0, 100)$)—an expected result since $\phi$ is close to one, as discussed in Section 2.

The second model is a newsvendor model with *two products* and limited order size. The first product has parameters $c = \eta = 300$, $p = \pi = 1500$, and demand is Uniform(100,400). The second product has parameters $c = \eta = 1000$, $p = \pi = 3000$, and demand is Uniform(50,150). With no limit on the order size, the problem decomposes as two independent problems with optimal order sizes respectively equal to 345.5 and 121.4; however, a limit of 300 on the total order size of the two products creates a constraint that ties the two problems, so an analytical solution is not available. We solve the problem by identifying an optimal scenario $(D_1^*, D_2^*)$ via a SAA version of (39)-(40) and solving the corresponding one-scenario problem. Again, we compare the result with a standard SAA solution of the two-stage formulation of this newsvendor model with two products, solved as a linear program. Table 3 displays the obtained optimal solutions and the corresponding CPU times for various sample sizes, ranging from $N = 100$ to $N = 50000$. Again, we observe that solving the problem via optimal scenarios yields exactly the same solution as the LP formulation of the two-stage problem; moreover, as discussed in Section 3.1 we see that the CPU times grow linearly when we solve the problem via (39)-(40), versus an exponential growth when we solve the two-stage problem as a linear program. While such a conclusion is based on a small experiment and as such cannot be generalized, the results suggest that it may be possible to develop an efficient algorithm that searches for optimal scenarios instead of searching for optimal solutions.

| Rep. | $(D^*, U^*)$ | $D^*/U^*$ | SAA solution of (18)-(23) | Analytical solution |
|---|---|---|---|---|
| 1 | (88.64, 0.41) | 218.02 | 218.02 | 214.73 |
| 2 | (172.9, 0.83) | 208.74 | 208.74 | 214.73 |
| 3 | (188.92, 0.89) | 211.42 | 211.42 | 214.73 |
| 4 | (190.0, 0.86) | 221.07 | 221.07 | 214.73 |
| 5 | (208.84, 0.96) | 217.30 | 217.30 | 214.73 |

Table 2: Optimal scenarios and optimal solutions for the newsvendor model with unreliable supplier.

| | solution via (39)-(40) | | | solution of the two-stage model | |
|---|---|---|---|---|---|
| $N$ | $(D_1^*, D_2^*)$ | $z^*$ | CPU time (sec) | $z^*$ | CPU time (sec) |
| 100 | (347.94, 101.98) | (198.02, 101.98) | 6.2 | (198.02, 101.98) | 0.008 |
| 1000 | (348.15, 98.42) | (201.58, 98.42) | 7.3 | (201.58, 98.42) | 0.12 |
| 5000 | (309.51, 99.28) | (200.72, 99.28) | 12.6 | (200.72, 99.28) | 1.8 |
| 10000 | (222.72, 99.38) | (200.62, 99.38) | 19.6 | (200.62, 99.38) | 6.2 |
| 50000 | (447.61, 98.51) | (201.49, 98.51) | 75.5 | (201.49, 98.51) | 153 |

Table 3: Optimal scenarios and optimal solutions for the newsvendor model with two products, multiple sample sizes.

## 7.2 Contextual Case: Synthetic Data Problems

### 7.2.1 Experiments Setting

To generate data for the problems described in Sections 7.2.2 and 7.2.3, we follow the approach in Kannan et al. (2022) (Section 4) to model the dependency among the clients' or locations' demands, $\xi = (\xi_j, j \in J) \in \mathbb{R}^J$, as a function of a covariate vector $x = (x_l, l \in L)$, for some finite index set $L$. Specifically, we assume that for all $j \in J$, the demand $\xi_j$ is given by

$$\xi_j = a_j + \sum_{l \in L} b_{j,l} \cdot (x_l)^p + \epsilon_j, \tag{90}$$

where $p > 0$ is a fixed *degree*, which we take in $\{0.5, 1, 2\}$, and $\epsilon_j \sim \mathcal{N}(0, \sigma_j^2)$ for all $j \in J$, with $\epsilon_j$ independent across $j \in J$. The parameters $(a_j, j \in J)$, $(b_{j,l}, (j,l) \in J \times L)$ and $(\sigma_j, j \in J)$ are consistent with those used in Section 4 of Kannan et al. (2022); see Appendix B for further details.

In our experiments, we consider instances with $|L| = 3$ covariates. The covariate vectors $x \in \mathbb{R}^L$ are modeled as i.i.d. draws from a multivariate folded (half-normal) distribution, following the procedure in Section 4 of Kannan et al. (2022). Specifically, we take $x_l = |\tilde{x}_l|$ for all $l \in L$, where the vector $\tilde{x} \in \mathbb{R}^L$ follows a multivariate normal distribution with zero mean and a covariance matrix sampled randomly from a Beta$(2, 2)$ distribution and subsequently rescaled to the interval $[-1, 1]$ (see Annex G of Kannan et al. (2022)).

We conduct experiments for values of the degree parameter $p \in \{0.5, 1, 2\}$ of the data generation procedure, and train the methods in Section 6 with $N$ samples of the pair $(x, \xi)$, where $N \in \{10^2, 10^3, 10^4\}$. Once training is complete, we evaluate the performance of each method by estimating its corresponding optimality gap using the procedure of Mak, Morton, and Wood (1999) (see Algorithm 1 in Kannan et al. 2022), summarized as follows:

1. Generate a sample of the covariate $x$.

2. Generate 1,000 samples of the demand $\xi$ conditional on the covariate value $x$, i.e., from the conditional distribution of $\xi|x$.

3. Compute the solution $z(x)$ of the method in Section 6 being tested, and compute the average cost of the solution $z(x)$ on the latter 1,000 samples.

4. Solve the problem $\min_{z \in Z} \mathbb{E}\left[G(z, \xi) \mid X = x\right]$ where the expected value is approximated using the 1,000 samples from the conditional distribution of $\xi|x$.

5. Compute the gap between the value obtained in step 3. minus the one from step 4.

6. Repeat 30 times steps 2. to 5., compute a confidence interval for the gap, and express it as a percentage of the value in step 4.

By using the above algorithm we obtain a normalized estimate of 99%-confidence upper bound for the optimality gap corresponding to the covariate $x$, which we denote $\widehat{B}_{99}(x)$. Since the data-driven solutions depend on the realization of the covariate sample $x$, we repeat 30 times this procedure —generating 30 covariates $x^1, \ldots, x^{30}$ and their corresponding upper bounds $\widehat{B}_{99}(x^1), \ldots, \widehat{B}_{99}(x^{30})$— and report the results using box plots of the latter 30 upper bounds.

We conducted our computational experiments by evaluating the performance of predict-then-optimize approaches, using CART and linear (LS) predictions, as well as application-driven approaches. Among the latter are the method described in Section 6.1.2, referred to as the M5 model (denoted M5+AD), and the method presented in Kannan et al. (2022), which corresponds to the policy in Equation (89), denoted as ER-SAA method. We also include the SAA- and KNN-based benchmark methods described in Section 6.2.2.

The hyper-parameters used in the decision tree models (CART, M5 and M5+AD) were the following. For both the shipment and resource allocation problems (Sections 7.2.2 and 7.2.3 below, respectively), the minimum number of samples per leaf was 25, and there was no restriction on the maximum depth of the tree. The value of the minimum number of samples was determined using k-fold cross-validation, a standard procedure that aims to maximize predictive performance while avoiding overfitting.

### 7.2.2 Problem 1: Two-Stage Resource Allocation

We first consider the two-stage resource allocation problem studied in Kannan et al. (2022). Let $I$ denote a set of resources and $J$ a set of clients. The first-stage decision involves determining the quantity $z_i \geq 0$ of resource $i$ to be ordered for each $i \in I$. The uncertain parameter is the demand $\xi_j$ of each client $j \in J$. Once demand is observed, two second-stage decisions are made: the amount $y^s_{i,j} \geq 0$ of resource $i$ allocated to client $j$, and the amount $y^c_j \geq 0$ of unmet demand for client $j$. The remaining parameters, assumed to be known, include the unit cost $c_i$ of resource $i$, the unit penalty cost $q_j$ for unmet demand of client $j$, the yield $\rho_i$ of resource $i$, and the service rate $\mu_{i,j}$ for allocating resource $i$ to client $j$. The problem can be formulated as follows:

$$\min_{z \in \mathbb{R}^I_+} \sum_{i \in I} c_i z_i + \mathbb{E}\left[Q(z, \xi)\right] \tag{91}$$

where the expectation is taken with respect to the random demand vector $\xi$, and the second-stage cost function $Q(z, \xi)$ is

$$Q(z, \xi) := \min_{y^s, y^c} \sum_{j \in J} q_j y^c_j \tag{92}$$

$$\text{s.t.} \quad \sum_{j \in J} y^s_{i,j} \leq \rho_i z_i \qquad \forall i \in I \tag{93}$$

$$\sum_{i \in I} \mu_{i,j} y^s_{i,j} + y^c_j \geq \xi_j \qquad \forall j \in J \tag{94}$$

$$y^s \in \mathbb{R}^{I \times J}_+, \ y^c \in \mathbb{R}^I_+ \tag{95}$$

In Figure 1 we show the results for the resource allocation problem in Section 7.2.2 with $|I| = 20$ resources, $|J| = 30$ clients and $|L| = 3$ covariates.

Our findings across the model degrees $p \in \{0.5, 1, 2\}$ are as follows. For all values of $p$, we see that the pure SAA approach, which ignores the contextual information, has poor performance. This is expected when there is some correlation between $\xi$ and $x$, as the knowledge of $x$ improves the estimation of $\xi$.

When the model degree is $p = 0.5$ or $p = 1$, methods AD, M5+AD, and ER-SAA show very good performance, even when the training data size is small. As more information becomes available, the kNN method also improves its performance. We can see that the predict-then-optimize policies CART and LS perform much worse than the other policies (except SAA), although the difference decreases as $N$ increases.

When the model degree is $p = 2$, we see a somewhat different behavior. The optimality gaps' variability increases due to the nonlinear form of the data. We also see that policies AD and ER-SAA end up performing worse (for larger $N$) than CART and similar to LS. The best performing methods are kNN and M5+AD, which take advantage of partitioning the data space.

Figure 5 in Appendix D shows a zoomed-in version of Figure 1 without CART, SAA and LS so as to emphasize the differences among the remaining four policies.

### 7.2.3 Problem 2: Two-Stage Shipment Planning

We also consider a two-stage shipment planning problem, originally introduced in Bertsimas and Kallus (2020). The model considers a set $I$ of warehouses and a set $J$ of demand locations. In the first stage, we determine the quantity $z_i \geq 0$ to be produced and stored at each warehouse $i \in I$, incurring a unit production cost $c > 0$. After the realization of demand $\xi_j$ at each location $j \in J$, two second-stage decisions are made: (i) the emergency production of $y^w_i$ units at warehouse $i$, at an elevated unit production cost $r > c$, and (ii) the quantity $y^s_{i,j} \geq 0$ of product shipped from warehouse $i$ to location $j$, incurring a unit shipment cost $s_{i,j}$.

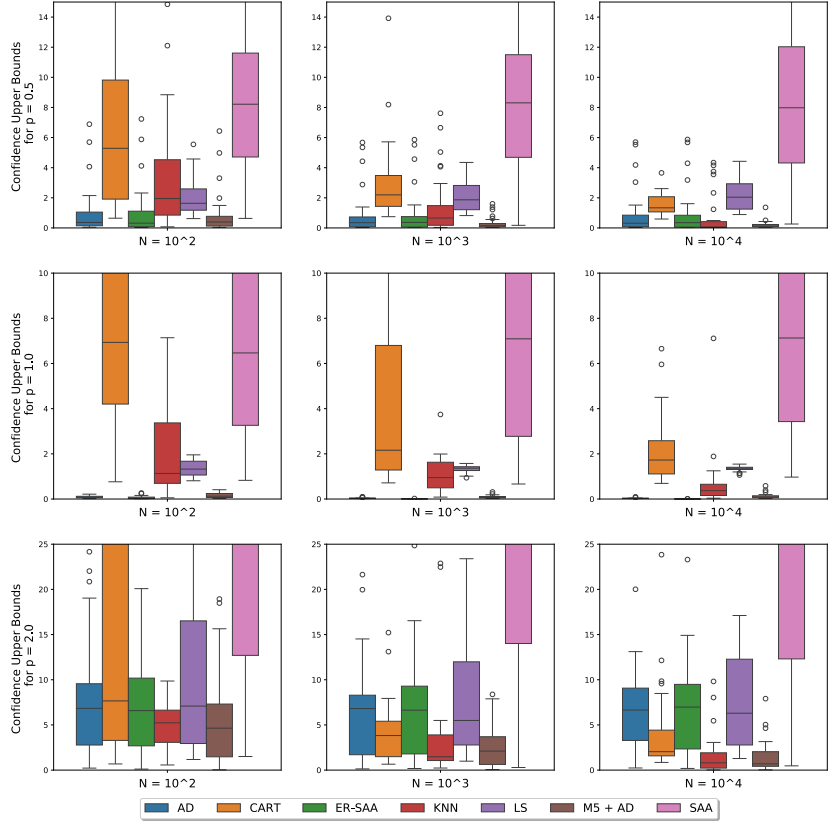

Figure 1: Comparison of methods proposed in Section 6 for the Problem 1 of resource allocation in Section 7.2.2. $p$ is the degree of the data generation procedure in (90), and $N$ is the number of samples with which each method is trained.

The problem is formulated as follows:

$$\min_{z \in \mathbb{R}_+^I} \quad c \sum_{i \in I} z_i + \mathbb{E}\left[Q(z, \xi)\right] \tag{96}$$

where the expectation is taken with respect to the random demand vector $\xi$, and the second-stage cost function $Q(z, \xi)$ is given by:

$$Q(z, \xi) := \min_{y^w, y^s} \quad r \sum_{i \in I} y_i^w + \sum_{i \in I} \sum_{j \in J} s_{i,j} y_{i,j}^s \tag{97}$$

$$\text{s.t.} \quad \sum_{i \in I} y_{i,j}^s \geq \xi_j \qquad \forall j \in J \tag{98}$$

$$\sum_{j \in J} y_{i,j}^s \leq z_i + y_i^w \qquad \forall i \in I \tag{99}$$

$$y^w \in \mathbb{R}_+^I, \ y^s \in \mathbb{R}_+^{I \times J} \tag{100}$$

In Figure 2 we show the results for the shipment planning problem in Section 7.2.3 with $|I| = 5$ warehouses,

$|J| = 12$ locations and $|L| = 3$ covariates. The conclusions are very similar to the previous problem, except that in this case the predict-then-optimize policies and the SAA policy perform much worse than AD policies and the conditional expectation ones even when $p = 2$. As before, Figure 6 in Appendix D shows a zoomed-in version of Figure 2 without CART, SAA and LS so as to emphasize the differences among the remaining four policies. We also see that, as in the previous problem, the methods that use an "optimalÂŽÂŽ pointwise forecast demonstrate in most cases similar or superior performance, compared to methods that consider multiple scenarios. In particular, the M5+AD method almost never loses to any of its competitors, regardless of the size of the training dataset, and has the best performance of all when $p = 2$. This supports our theoretical results, showing that a two-stage stochastic program can actually be solved with only one scenario.

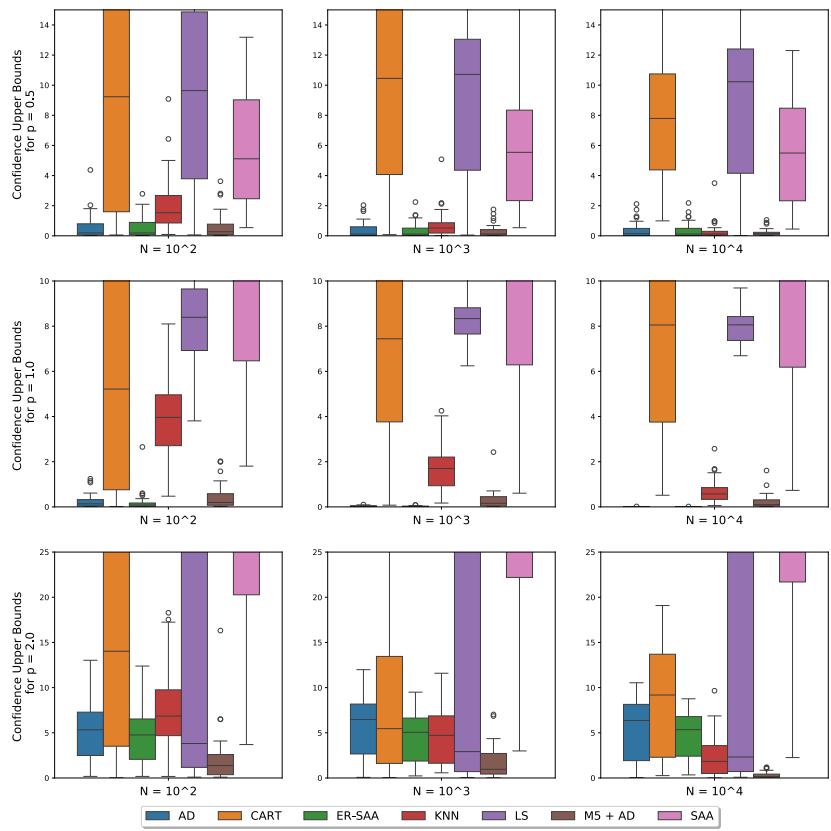

Figure 2: Comparison of methods proposed in Section 6 for Problem 2 of shipment planning in Section 7.2.3. $p$ is the degree of the data generation procedure in (90), and $N$ is the number of samples with which each method is trained.

## 7.3 Real-World Data Problem: Bike Sharing Reallocation

To evaluate our methodology on a realistic setting with real-world data, we consider the bike reallocation problem studied in Cavagnini (2019), which uses data from the San Francisco, CA, bike-sharing system. The data are publicly available through the open-source repository `https://www.kaggle.com/datasets/`

`benhamner/sf-bay-area-bike-share/data` and describe a system with 350 bicycles and 34 stations (see the station distribution in Figure 3). The repository comprises four datasets containing multiple fields, including station status, weather conditions (temperature, wind speed, and humidity), and trip information, covering the period from August 2013 to August 2015.

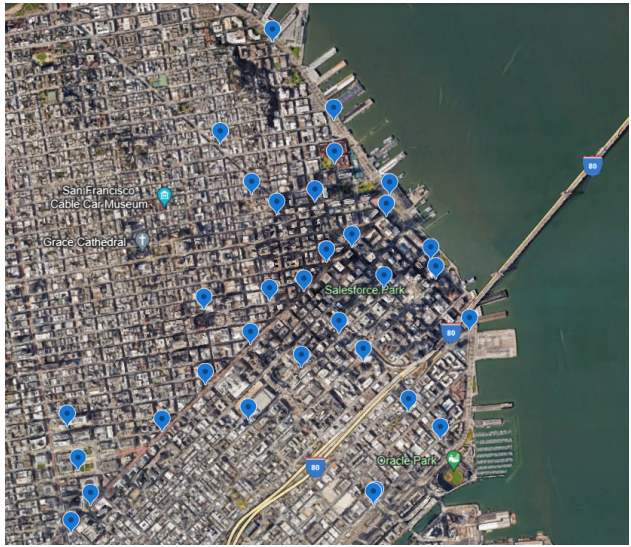

Figure 3: San Francisco's stations of bike sharing system

Cavagnini (2019) propose a two-stage stochastic model to address the bike-sharing problem in San Francisco. Although their formulation does not incorporate features or covariate information, it nonetheless fits within the broader framework of data-driven stochastic programming with covariates considered in this paper. Their model includes integer decision variables; however, in order to apply our methodology, we relax these integrality constraints and restrict attention to continuous variables.

### 7.3.1 Problem Formulation

With a slight abuse of notation, let $\mathcal{I} = \{1, 2, \ldots, |\mathcal{I}|\}$ denote the set of stations, where station $|\mathcal{I}|$ corresponds to the bike depot. The problem begins at the start of the day, when the decision-maker must determine the number of bikes $x_i$ to deliver from the depot to each station, incurring a unit delivery cost of $f_i$. This decision occurs before the demand $\xi_i$ at each station $i \in \mathcal{I} \setminus \{|\mathcal{I}|\}$ is known. In our setting, we assume that the decision-maker can forecast station-level demand using covariate information, which will be described later. Each station $i$ has a maximum capacity of $Q_i$ bikes and requires a minimum allocation of $\underline{x}_i$ bikes at the beginning of the day to guarantee service feasibility. Once the demand $\xi_i$ is realized, the service provider performs a rebalancing operation at the end of the day, redistributing $y_{i,i+1}$ bikes from station $i \in \mathcal{I} \setminus \{|\mathcal{I}|\}$ to the next station along a fixed route-which, without loss of generality, we assume to be station $i+1$-at a relocation cost of $t_{i,i+1}$.

We also assume that the route begins and ends at the bike depot. The reallocation of bikes is carried out using a truck with total capacity $C$. The objective of the bike-sharing provider is to avoid situations in which a user attempts to return a bike to a station that is already full or to rent a bike from an empty station. At the same time, the decision-maker seeks to minimize the number of bikes relocated, thereby reducing the risk of bike damage. To this end, Cavagnini (2019) introduce *starvation* and *congestion* terms.

Starvation is captured by the variable $I_i^-$, representing stock-outs at station $i \in \mathcal{I} \setminus \{|\mathcal{I}|\}$, and penalized through the *stock-out penalty* $p_i$. Congestion is represented by two additional terms for each station $i \in$

| *Sets* | |
|---|---|
| $\mathcal{I} = \{1, \ldots, |\mathcal{I}|\}$ | Set of stations, where the depot corresponds to station $|\mathcal{I}|$ |
| *Parameters* | |
| $\underline{x}_i$ | Minimum number of bikes to allocate at station $i \in \mathcal{I} \setminus \{|\mathcal{I}|\}$ |
| $\bar{I}_{|\mathcal{I}|0}$ | Depot capacity |
| $Q_i$ | Capacity of station $i \in \mathcal{I} \setminus \{|\mathcal{I}|\}$ |
| $\bar{I}_{i0}$ | Initial number of bikes at station $i \in \mathcal{I} \setminus \{|\mathcal{I}|\}$ |
| $C$ | Capacity of the relocation truck |
| $p_i$ | Stock-out penalty at station $i \in \mathcal{I} \setminus \{|\mathcal{I}|\}$ |
| $c_i$ | Excess penalty at station $i \in \mathcal{I} \setminus \{|\mathcal{I}|\}$ |
| $\frac{c_i}{Q_i}$ | Penalty for each extra bike placed at station $i \in \mathcal{I} \setminus \{|\mathcal{I}|\}$ after rebalancing |
| $f_i$ | Allocation cost at station $i \in \mathcal{I} \setminus \{|\mathcal{I}|\}$ |
| $t_{i,i+1}$ | Rebalancing cost of moving bikes from station $i$ to $i+1$, $i \in \mathcal{I} \setminus \{|\mathcal{I}|\}$ |
| $\xi_i$ | Net demand at station $i \in \mathcal{I} \setminus \{|\mathcal{I}|\}$ |
| *Variables* | |
| $x_i$ | Bikes allocated to station $i \in \mathcal{I} \setminus \{|\mathcal{I}|\}$ in the first stage |
| $y_{i,i+1}$ | Bikes distributed from station $i$ to $i+1$ in the second stage, $i \in \mathcal{I} \setminus \{|\mathcal{I}|\}$ |
| $I_i$ | Inventory of bikes at station $i \in \mathcal{I} \setminus \{|\mathcal{I}|\}$ |
| $I_i^+$ | Surplus bikes at station $i \in \mathcal{I} \setminus \{|\mathcal{I}|\}$ |
| $I_i^-$ | Stock-out bikes at station $i \in \mathcal{I} \setminus \{|\mathcal{I}|\}$ |
| $B_i$ | Extra bikes at station $i \in \mathcal{I} \setminus \{|\mathcal{I}|\}$ during rebalancing |
| $B_i^+$ | Balance of extra bikes at station $i \in \mathcal{I} \setminus \{|\mathcal{I}|\}$ during rebalancing |
| $E_i$ | Excess inventory of bikes at station $i \in \mathcal{I} \setminus \{|\mathcal{I}|\}$ during rebalancing |
| $E_i^+$ | Excess bikes at station $i \in \mathcal{I} \setminus \{|\mathcal{I}|\}$ during rebalancing |

Table 4: Sets, parameters, and variables for Problem 3 in Section 7.3.

$\mathcal{I} \setminus \{|\mathcal{I}|\}$: the *extra inventory* variable $B_i^+ \geq 0$, which measures the number of bikes exceeding the initially allocated amount, and the *excess inventory* variable $E_i^+ \geq 0$, which measures the number of bikes surpassing the station capacity $Q_i$. Nonzero values of these variables are penalized with unit costs of $c_i$ and $c_i/Q_i$, respectively

With this, the problem formulation is as follows:

$$\min_x \quad 3\sum_{i \in \mathcal{I}} f_i x_i + \mathbb{E}\left[Q(x, \xi)\right] \tag{101}$$

$$\text{s.t.} \quad x_i \geq \underline{x}_i \qquad \forall i \in \mathcal{I} \setminus \{|\mathcal{I}|\} \tag{102}$$

$$\bar{I}_{i0} + x_i \leq Q_i \qquad \forall i \in \mathcal{I} \setminus \{|\mathcal{I}|\} \tag{103}$$

$$\sum_{i \in \mathcal{I} \setminus \{|\mathcal{I}|\}} x_i \leq \bar{I}_{|\mathcal{I}|0} \tag{104}$$

$$x_i \geq 0 \qquad \forall i \in \mathcal{I} \tag{105}$$

where the second-stage cost function $Q(x,\xi)$ is

$$Q(x,\xi) = \min \sum_{i \in \mathcal{I} \setminus \{|\mathcal{I}|\}} \left( t_{i,i+1} y_{i,i+1} + \frac{c_i}{Q_i} B_i^+ + c_i E_i^+ + p_i(-I_i^-) \right) \quad (106)$$

$$\text{s.t.} \quad y_{i,i+1} \leq C \quad \forall i \in \mathcal{I} \setminus \{|\mathcal{I}|\} \quad (107)$$

$$I_{|\mathcal{I}|} = \bar{I}_{|\mathcal{I}|0} - \sum_{i \in \mathcal{I} \setminus \{|\mathcal{I}|\}} x_i + y_{|\mathcal{I}|-1,|\mathcal{I}|} \quad (108)$$

$$I_{|\mathcal{I}|} \leq \bar{I}_{|\mathcal{I}|0} \quad (109)$$

$$I_1 = \bar{I}_{|\mathcal{I}|0} + x_1 - \xi_1 - y_{1,2} \quad (110)$$

$$I_i = \bar{I}_{i0} + x_i - \xi_i + y_{i-1,i} - y_{i,i+1} \quad \forall i \in \mathcal{I} \setminus \{1, |\mathcal{I}|\} \quad (111)$$

$$I_i^+ = \max\{0, I_i\} \quad \forall i \in \mathcal{I} \setminus \{|\mathcal{I}|\} \quad (112)$$

$$I_i^- = \min\{0, I_i\} \quad \forall i \in \mathcal{I} \setminus \{|\mathcal{I}|\} \quad (113)$$

$$E_i = I_i^+ - Q_i \quad \forall i \in \mathcal{I} \setminus \{|\mathcal{I}|\} \quad (114)$$

$$E_i^+ = \max\{0, E_i\} \quad \forall i \in \mathcal{I} \setminus \{|\mathcal{I}|\} \quad (115)$$

$$B_i = I_i^+ - x_i - \bar{I}_{i0} - E_i^+ \quad \forall i \in \mathcal{I} \setminus \{|\mathcal{I}|\} \quad (116)$$

$$B_i^+ = \max\{0, B_i\} \quad \forall i \in \mathcal{I} \setminus \{|\mathcal{I}|\} \quad (117)$$

$$I_i, B_i, E_i \in \mathbb{R} \quad \forall i \in \mathcal{I} \setminus \{|\mathcal{I}|\} \quad (118)$$

$$y_{i,i+1}, I_i^+, B_i^+, E_i^+ \geq 0 \quad \forall i \in \mathcal{I} \setminus \{|\mathcal{I}|\} \quad (119)$$

$$I_i^- \leq 0 \quad \forall i \in \mathcal{I} \setminus \{|\mathcal{I}|\} \quad (120)$$

Constraints (112), (113), (115), and (117) can be readily linearized by introducing auxiliary decision variables to represent the max and min terms.

We now describe the components of the above model, following closely the discussion in Cavagnini (2019). The objective functions (101) and (106) minimize the total expected penalty, obtained by summing delivery, rebalancing, extra and excess inventory, and stock-out costs. Constraints (102) require that the number of bikes delivered to each station be at least the specified minimum. Constraints (103) ensure that the sum of allocated bikes and initial availability at each station does not exceed station capacity. Constraints (104) impose that the total number of allocated bikes across all stations is no greater than the available quantity at the depot. Constraints (107) guarantee that the number of bikes transported by the truck during rebalancing never exceeds its capacity. Constraints (108) enforce that the depot inventory at the end of each scenario equals the initial availability plus the amount received from the last visited station, minus the amounts delivered to stations. Similarly, constraints (109) restrict the depot inventory at the end of the rebalancing period to not exceed depot capacity.

The flow balance constraints differ between the first station and subsequent stations on the route. Constraints (110) ensure that, for the first station, the final inventory equals the initial availability plus the amount received from the depot, minus the demand served and the number of bikes redistributed downstream. Constraints (111) define the inventory at any other station as a function of its initial level, the number allocated, the bikes withdrawn or returned, and the amount redistributed further along the route. Constraints (112) and (113) determine the surplus and stock-out quantities at each station. When a valet service is available (allowing bikes to be returned even to full stations) Constraints (114) and (115) compute the number of bikes exceeding station capacity. Finally, Constraints (116) and (117) identify cases in which the number of bikes after rebalancing exceeds the initial allocation but remains within station capacity.

### 7.3.2 Results

We now present the results of applying the methods described in Section 6 to the two-stage problem introduced above. The uncertain parameter is the daily bike demand at each station, for which our methodology uses additional covariate information to generate forecasts. In our experiments, we use six features: a binary indicator denoting business versus non-business days, the week number, the month of the year, the year, the level of precipitation (in inches), and the wind direction. After aggregating demand and features on a daily basis, the dataset contains approximately 730 observations for prediction and optimization.

The hyper-parameters used in the decision tree models (CART, M5 and M5+AD) were the following: the minimum number of samples per leaf was 1 and the maximum tree depth was set to 8. The value of the minimum number of samples was determined using k-fold cross-validation, a standard procedure that aims to maximize predictive performance while avoiding overfitting. Also, the maximum tree depth of 8 was determined based on the fact that there were only 730 data points, and a tree with depth $D$ can generate up to $2^D$ leaves.

To estimate the optimality gap of each method we compute its out-of-sample performance by randomly splitting the dataset in a training and test set, keeping a standard proportion of 80:20 of training and test size, and we take the *out-of-sample cost* as the average of the bike-sharing problem cost over the testing data. Since results may depend on the specific realization of the training and test partitions, we repeat this procedure over 10 random splits, preserving the same proportions.

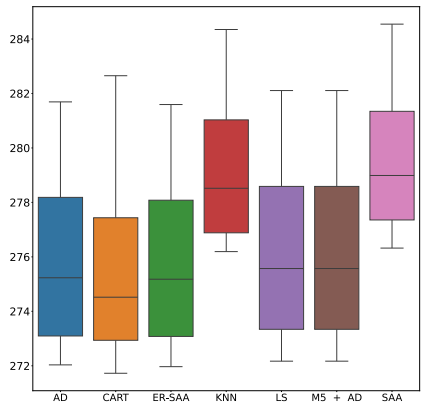

Figure 4: Average out-of-sample costs for 10 different train and test sets for methods proposed in Section 6.

Figure 4 reports the out-of-sample costs for all methods under the multiple-replications procedure. We observe that all methods, except KNN and SAA, perform similarly, while KNN and SAA exhibit clearly inferior results. The similarity among the five best-performing methods may be explained by the relatively small dataset (730 observations), as some of these methods require larger sample sizes to realize their full potential. Consistent with previous experiments, the weaker performance of SAA can be attributed to its inability to exploit contextual information. Overall, this example confirms that methods based on well-chosen pointwise forecasts, such as AD and M5+AD, achieve comparable or superior performance relative to approaches that rely on multiple scenarios.

**Remark:** Cavagnini (2019) proposes the two-stage bike relocation model presented in (101)-(120) that assumes static demand, but in the published version of that work (Cavagnini, Maggioni, Bertazzi, and Hewitt, 2024) the authors attempt to address that limitation by proposing a small model variation that requires estimating additional random quantities—namely, the maximum number of consecutive bikes that are returned to a station before a bike is withdrawn, and the maximum number of consecutive bikes that are withdrawn from a station before a bike is returned. Unfortunately, the limited size of the dataset does not allow for an accurate estimation of those quantities (68, i.e., twice the number of stations as we now have two new variables per station); in fact, Cavagnini et al. (2024) generate random samples from the empirical distribution. As our goal in this example is to truly estimate the random quantities directly from the data—by using contextual information—one possible way to circumvent the problem is to define such quantities as fractions of the demand in each node. While such an approach bypasses the need for more data—since it only requires estimating the demand, as in the original model—it leads to functional dependencies among the right-hand side of different constraints, which as discussed in Section 5.3 is not respected by the one-scenario approach. In Appendix E we present the results, where we see that indeed the AD and M5+AD approaches perform worse than other methods. In fact, in that case SAA is among the best performers, thereby suggesting that the contextual information is of little use in such a model. In summary, to properly solve the model in Cavagnini et al. (2024) using contextual information it would be necessary to have a larger dataset that would allow us to estimate all required quantities directly from the data.

# 8 Conclusions

Stochastic optimization problems are defined in terms of the (possibly unknown) distributions of the underlying random variables. Accordingly, methods to solve such problems typically estimate the input distributions and then apply some some scenario generation/reduction technique, perhaps combined with an approach that allows for the decomposition of the problem across scenarios. In this paper we have considered an alternative application-driven approach whereby only *pointwise estimates* are required, when the problem to be solved belongs to the class of two-stage stochastic programs with fixed recourse and fixed costs. The basis for the proposed approach is a novel result that shows that, for that class of problems, it suffices to use *one* scenario, in the sense that solving the problem with that single scenario yields the same solution as the original problem.

In the setting of optimization problems with contextual information, the task of estimating input distributions becomes more difficult as it requires the estimation of conditional distributions for any given value of the contextual information. On the other hand, *end-to-end learning* techniques proposed in the literature have proven very valuable in the contextual information setting by combining the estimation and optimization steps. By using our "optimal scenario" result, our integrated learning and optimization method uses problem information to determine pointwise forecasts that provide the best parametric approximation of that (unknown) scenario.

The main goal of this work is to show that our novel pointwise approach provides a practical alternative way to solve two-stage problems with contextual information, which by-passes the need for estimating conditional distributions. Our numerical results corroborate that idea. There is of course much room for improvements, especially regarding the development of specialized techniques to solve the bilevel models that are part of the method, and the use of other machine learning methods such as neural networks within our setting. We hope our work will stimulate further research on these topics. We also hope that our one-scenario result can spur new research on alternative methods for stochastic optimization problems (with or without contextual information) that search over the space of scenarios rather than over the space of solutions.

# Appendices

## A  Proof of expression (18)

Equation (17) can be written as

$$\int_0^1 uF(zu)dG(u) = \phi\mathbb{E}[U], \tag{121}$$

where $F$ is the distribution of the demand $D$ and $G$ is the distribution of the reliability factor $U$. Consider the case where

$$U \sim \text{Uniform}(0,1)$$
$$D \sim \text{Uniform}(0,b).$$

Then, we have

$$F(y) = \begin{cases} 0 & y \leq 0 \\ \frac{y}{b} & 0 < y < b \\ 1 & y \geq b \end{cases}$$

and thus the left-hand side of (121) is

$$\int_0^1 uF(zu)dG(u) = \int_0^1 u\mathbf{1}_{\{0<uz\leq b\}}\frac{zu}{b}\,du + \int_0^1 u\mathbf{1}_{\{zu>b\}}\,du$$
$$= \frac{1}{b}\int_0^{\min(b/z,1)} zu^2\,du + \int_{\min(b/z,1)}^1 u\,du$$
$$= \frac{1}{b}\left(\frac{zu^3}{3}\right)\Big|_0^{\min(b/z,1)} + \frac{u^2}{2}\Big|_{\min(b/z,1)}^1. \tag{122}$$

Note that for $z \leq 0$ the left hand side of (121) is equal to zero so that equation has no solution. Suppose now that $0 < z \leq b$. Then we have $\min(b/z,1) = 1$ and thus from (122) it follows that

$$\int_0^1 uF(zu)dG(u) = \frac{z}{3b}.$$

Thus, in (121) we have

$$\frac{z}{3b} = \frac{\phi}{2}$$

i.e., $z = \frac{3}{2}\phi b$. Notice that the condition $z \leq b$ is satisfied as long as $\phi \leq 2/3$.

Suppose now that $z > b$. Then we have $\min(b/z,1) = b/z$ and thus from (122) it follows that

$$\int_0^1 uF(zu)dG(u) = \frac{b^2}{3z^2} + \left(\frac{1}{2} - \frac{b^2}{2z^2}\right).$$
$$= \frac{1}{2} - \frac{b^2}{6z^2}.$$

Thus, in (121) we have

$$\frac{1}{2} - \frac{b^2}{6z^2} = \frac{\phi}{2}$$

that is

$$\frac{b^2}{3z^2} = 1 - \phi$$

and thus the solution is
$$z = \frac{b}{\sqrt{3(1-\phi)}}$$

Notice that the condition $z > b$ is satisfied as long as $3(1-\phi) < 1$, i.e., $\phi > 2/3$.

# B  Data generation details

The parameters $(a_j, j \in J)$, $(b_{j,l}, (j,l) \in J \times L)$ and $(\sigma_j, j \in J)$ in the data generation procedure (90) are taken in Kannan et al. (2022) as follows:

$$a_j = 50 + 5\delta_{j,0} \tag{123}$$
$$b_{j,1} = 10 + \delta_{j,1} \tag{124}$$
$$bj,2 = 5 + \delta_{j,2} \tag{125}$$
$$bj,3 = 2 + \delta_{j,2} \tag{126}$$

where $\{\delta_{j,0}\}_{j \in J}$ are i.i.d. samples from the standard normal distribution, and $\{\delta_{j,1}\}_{j \in J}$, $\{\delta_{j,2}\}_{j \in J}$, $\{\delta_{j,3}\}_{j \in J}$ are i.i.d. samples from a uniform distribution $U(-4, 4)$. Also, $\sigma_j = \sigma = 5$ for all $j$ in $J$.

# C  Estimation of the stations' daily demand

To calculate the net demand for each station $d_i$, we need to forecast the daily demand for the entire system. This is based on the features $X$ that provide information on the daily weather in San Francisco. Subsequently, an adjacency matrix $P_{ij}$ is computed to estimate the proportion of daily trips between stations. The demand for each station is determined by subtracting the number of bikes returned $d_i^r$ from the number of bikes withdrawn $d_i^w$ at each station, following the procedure proposed by Cavaginini. To estimate the daily demand for stations using different methods, we first calculate the overall demand using the methods outlined in Section 6. Then, we calculate the estimated number of bikes withdrawn and returned to each station by multiplying the adjacency matrix by the forecast obtained from the methods proposed. Finally, the net demand is obtained by subtracting the estimated number of withdrawn and returned bikes (i.e. $d_i = d_i^w - d_i^r$). The data that supports the experiments conducted within problem in Section 7 and Appendix E is available at: https://www.kaggle.com/datasets/benhamner/sf-bay-area-bike-share.

**Data treatment and forecasting.** The dataset provides four datasets where information on the trips, theater, and stations' status and information are given. The trips dataset contains information on all the trips within stations, indicating the station and time of beginning and end of each trip. The trip data was grouped to present daily trips for the whole system. The weather dataset contains daily information about weather conditions such as temperature, wind speed, and others. The dataset comprises continuous and categorical variables, however, a vectorization of categorical variables is proposed to address this issue. A binary variable is created to indicate whether a particular day is a holiday. Using station data, a variable has been created to indicate the available docks within stations for each day. Lastly, the station data is used to estimate optimization model parameters such as the capacity of each station. In regards to forecasting, after conducting feature engineering and a feature importance procedure, we selected 6 features to predict the demand for the methods outlined in Section 6. Based on these features, the forecast aims to predict the number of bikes in the system.

**Adjacent matrix calculation and stations' estimated net demand.** We used the trip data to compute the adjacency matrix $P_{ij}$. This involved counting the number of bikes that departed from each station $i$ and ended at station $j$, which is represented in the $i,j$ position of the matrix obtaining a matrix $C_{ij}$. It is important to mention that the diagonal of the matrix represents the trips that start and end in the same station. To obtain the proportion of the trips we divide the matrix $C_{ij}$ with the total amount of trips of the system in the whole dataset ($P_{ij} := C_{ij}/\sum_{k,l} C_{kl}$). Finally, given the proportion of trips and the daily system forecast for each method, $\hat{y}^t_{sys}$, we obtain the estimated trips within a day of operation by simply multiplying the estimated daily demand by the adjacency matrix ($\hat{y}^t_{sys} P_{ij}$). Finally, given the estimate of the withdrawn and returned bikes for each day and each method, we simply compute the estimated net demand $d_i$ by subtracting the estimated number of withdrawn and returned (i.e. $d_i^w - d_i^r$), which are computed by summing the rows and columns of each station respectively, obtaining the net demand for each station $d_i$.

## D  Zoomed-in versions of Figures 1 and 2

We present below zoomed-in versions of Figures 1 and 2 that display only the best four contenders, i.e. we remove the SAA, LS and CART policies.

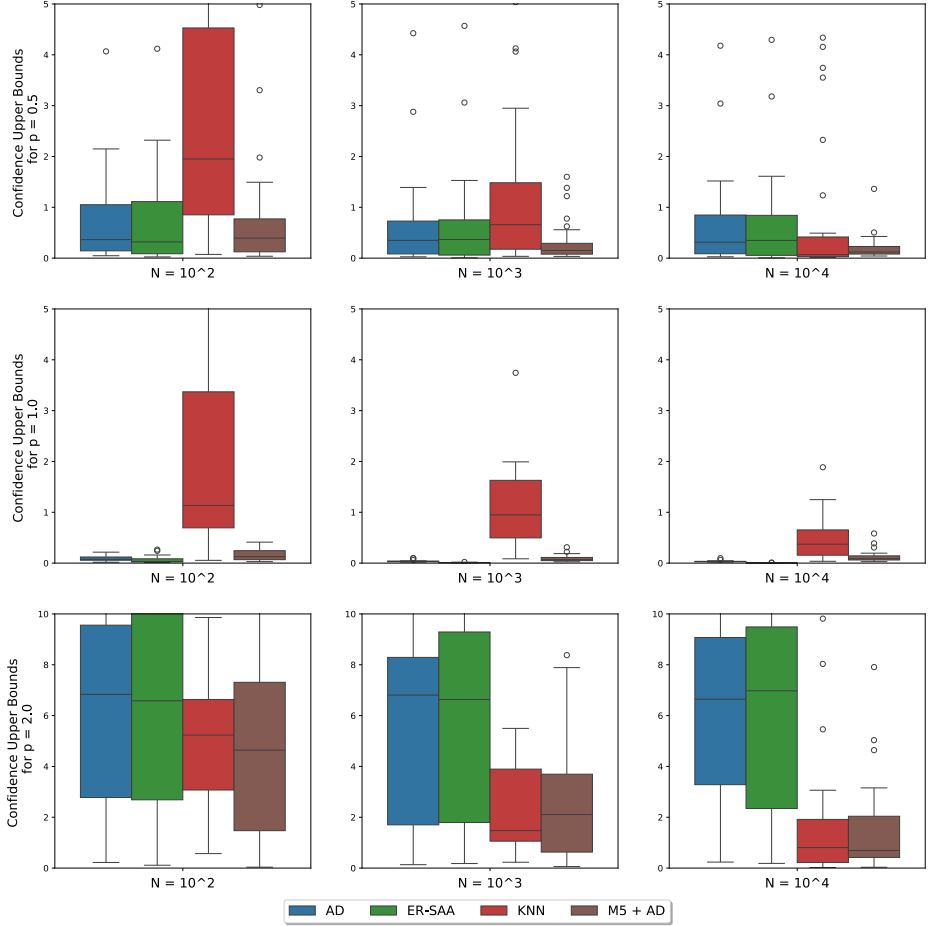

Figure 5: Detailed comparison of methods proposed in Section 6 for the Problem 1 of resource allocation in Section 7.2.2. $p$ is the degree of the data generation procedure in (90), and $N$ is the number of samples with which each method is trained.

# E  A reformulation of the bike sharing rellocation model

Cavagnini (2019) proposed the two-stage bike relocation model presented in (101)-(120) that assumed static demand, but Cavagnini et al. (2024) attempted to address this limitation by proposing a small model variation. The proposed model ensures that there are always at least as many bikes available as the maximum number of consecutive bikes that have been withdrawn from a station before being returned. It also ensures that there are always at least as many free docks available as the maximum number of consecutive bikes that have been returned to a station before being withdrawn. This approach helps to prevent congestion and shortages, especially during times of peak consecutive bike withdrawals and returns. To archive that, Cavagnini et al. (2024) introduced extra stochastic parameters $g_i$ and $h_i$ representing the maximum number of consecutive bikes withdrawn from station $i \in \mathcal{I} \setminus \{|\mathcal{I}|\}$ before a return occurs, and the maximum number of consecutive bikes returned station $i \in \mathcal{I} \setminus \{|\mathcal{I}|\}$ before a withdrawal occurs. These parameters are obtained by using Montecarlo sampling using the historical demand, however, in our setting, we considered them as fixed parameters due to feasibility, which depends on the estimate of the total number of bikes returned $d_i^r$ and the estimate of the total of bikes withdrawn $d_i^w$ at each station.

**Estimating the new parameters.** To determine the parameters, we calculated the total number of bikes withdrawn, $W_{i,t}$, and returned, $R_{i,t}$, as well as the maximum number of consecutive bikes withdrawn, $w_{i,t}$, and returned, $r_{i,t}$, at each station $i \in \mathcal{I} \setminus \{|\mathcal{I}|\}$ on each day $t \in T$ within the trip dataset. With this data, we calculated $\alpha_{i,t} := w_{i,t}/W_{i,t}$ and $\beta_{i,t} := r_{i,t}/R_{i,t} \quad \forall_{i \in \mathcal{I} \setminus \{|\mathcal{I}|\}, t \in T}$, representing the proportion of the maximum number of consecutive bikes withdrawn and returned over the total withdrawn and returned bikes at each station. Then, we calculate the vectors $\bar{\alpha}_i := \frac{1}{T}\sum_{t \in T} \alpha_{i,t}$ and $\bar{\beta}_i := \frac{1}{T}\sum_{t \in T} \beta i, t$ that represent the expected proportion of maximum number of consecutive bikes withdrawn and returned over the total number of bikes withdrawn and returned. Finally, we estimate $h_i$ and $g_i$ as a proportion of the estimated withdrawn and returned demand, that is, $h_i := \bar{\alpha}_i d_i^w$ and $g_i := \bar{\beta}_i d_i^r$ where $d_i^r$ and $d_i^w$ are obtained by using the forecast methods of Section 6 and the methodology detailed in Appendix C.

**Problem formulation.** The model proposed by Cavagnini et al. (2024) is similar to the proposed in Cavagnini (2019), however, they introduce some parameters and constraints to address the limitation of the static demand. The problem considers a set $\mathcal{I} = \{1, 2, \ldots, |\mathcal{I}|\}$ of stations with capacity of $Q_i$ bikes, where station $|\mathcal{I}|$ is the bike depot. The bike planning starts at the beginning of the day, where, before knowing the demand $\xi_i$ of bikes at each station $i \in \mathcal{I} \setminus \{|\mathcal{I}|\}$, we have to allocate a number $x_i$ of bikes to deliver from the depot to the station having a fixed cost of $f_i$ per bike.

Later in the day, the demand for each station $\xi_i$ is realized. At the end of the day, the service provider rebalances the bikes, having to re-distribute $y_{i,i+1}$ bikes from station $i \in \mathcal{I} \setminus \{|\mathcal{I}|\}$ to the next station $i+1$ on the fixed route, incurring a moving cost of $t_{i,i+1}$. The bikes are redistributed using a truck with a total capacity of $C$ bikes. The bike provider aims to avoid situations where a user wants to return a bike to a full station or needs to rent a bike from an empty station. The decision maker also wants to minimize the number of bikes redistributed. To achieve this, Cavagnini suggests using starvation and a congestion term. The starvation is parametized by the variable $I_i^-$ and a *stock-out penalty* $p_i$, for each station $i \in \mathcal{I} \setminus \{|\mathcal{I}|\}$. The congestion is measured with two terms at each station $i$ in $\mathcal{I} \setminus \{|\mathcal{I}|\}$: the *extra inventory* term $B_i^+ \geq 0$ measures the number of bikes beyond the number initially allocated, and the *excess inventory* term $E_i^+ \geq 0$ measures the number of bikes over the station's capacity $Q_i$. Any non-negative value of these is penalized with a unit cost of $c_i$ and $c_i/Q_i$, respectively. Finally, the model introduces the stochastic parameters $g_i$ and $h_i$ are used to prompt the model to determine a target inventory level that is higher than the maximum number of consecutive withdrawn bikes, $g_i$, and to ensure that there are at least $h_i$ free docks. To ensure that we can find a feasible solution even if the total number of bikes withdrawn and returned consecutively exceeds the station capacity, we introduce the variables $a_i$ to represent the difference between the number of allocated bikes $x_i$ and $g_i$ at station $i \in \mathcal{I} \setminus \{|\mathcal{I}|\}$. Similarly, the variables $b_i$ represent the difference between the number of free racks $Q_i - x_i$ and $h_i$ at station $i \in \mathcal{I} \setminus \{|\mathcal{I}|\}$.

The problem formulation is as follows:

| Sets | |
|---|---|
| $\mathcal{I} = \{1, \ldots, |\mathcal{I}|\}$ | Set of stations, where the depot is station $|\mathcal{I}|$ |
| *Parameters* | |
| $\bar{I}_{|\mathcal{I}|0}$ | Depot capacity |
| $Q_i$ | Maximum capacity of station $i \in \mathcal{I} \setminus \{|\mathcal{I}|\}$ |
| $\bar{I}_{i0}$ | Initial availability of bikes at station $i \in \mathcal{I} \setminus \{|\mathcal{I}|$ |
| $C$ | Maximum capacity of the relocation truck |
| $p_i$ | bikes' stock-out penalty at station $i \in \mathcal{I} \setminus \{|\mathcal{I}|\}$ |
| $c_i$ | Excess penalty at station $i \in \mathcal{I} \setminus \{|\mathcal{I}|\}$ |
| $\frac{c_i}{Q_i}$ | Penalty associated to extra bikes placed at station $i \in \mathcal{I} \setminus \{|\mathcal{I}|$ after the re-balancing period |
| $f_i$ | Allocation cost at station $i \in \mathcal{I} \setminus \{|\mathcal{I}|\}$ |
| $t_{i,i+1}$ | Re-balancing cost to allocate bikes from station $i$ to $i+1$, $i \in \mathcal{I} \setminus \{|\mathcal{I}|\}$ |
| $\xi_i$ | Net demand of station $i \in \mathcal{I} \setminus \{|\mathcal{I}|\}$ |
| $g_i$ | Maximum number of consecutive bikes withdrawn from station $i \in \mathcal{I} \setminus \{|\mathcal{I}|\}$ |
| $h_i$ | Maximum number of consecutive bikes returned from station $i \in \mathcal{I} \setminus \{|\mathcal{I}|\}$ |
| *Variables* | |
| $x_i$ | Amount of bikes to allocate at station $i \in \mathcal{I} \setminus \{|\mathcal{I}|\}$ in the first stage |
| $a_i$ | Slack units between the number of allocated bikes in the first stage and the maximum number of consecutive bikes withdrawn at station $i \in \mathcal{I} \setminus \{|\mathcal{I}|\}$ |
| $b_i$ | Slack units between the number of available racks and the maximum consecutive bikes returned at station $i \in \mathcal{I} \setminus \{|\mathcal{I}|\}$ |
| $y_{i,i+1}$ | Amount of bikes to distribute from station $i$ to $i+1$ in the second stage, $i \in \mathcal{I} \setminus \{|\mathcal{I}|\}$ |
| $I_i$ | Inventory or balance of bikes at station $i \in \mathcal{I} \setminus \{|\mathcal{I}|\}$ |
| $I_i^+$ | Surplus of bikes at station $i \in \mathcal{I} \setminus \{|\mathcal{I}|\}$ |
| $I_i^-$ | Number of stock-out bikes at station $i \in \mathcal{I} \setminus \{|\mathcal{I}|\}$ |
| $B_i$ | Extra bikes balance at station $i \in \mathcal{I} \setminus \{|\mathcal{I}|\}$ during the re-balancing period |
| $B_i^+$ | Balance of of extra bikes at station $i \in \mathcal{I} \setminus \{|\mathcal{I}|\}$ during the re-balancing period |
| $E_i$ | Excess inventory balance of bikes at station $i \in \mathcal{I} \setminus \{|\mathcal{I}|\}$ during the re-balancing period |
| $E_i^+$ | Excess of bikes at station $i \in \mathcal{I} \setminus \{|\mathcal{I}|\}$ during the re-balancing period |

Table 5: Variables, parameters, and sets for Problem 4 in AAAAA.

$$\min_x \quad \sum_{i \in \mathcal{I}} f_i x_i + \mathbb{E}\left[Q(x, \xi)\right] \tag{127}$$

$$\text{s.t.} \quad x_i \geq g_i - a_i \qquad \forall i \in \mathcal{I} \setminus \{|\mathcal{I}|\} \tag{128}$$

$$Q_i - x_i \geq h_i - b_i \qquad \forall i \in \mathcal{I} \setminus \{|\mathcal{I}|\} \tag{129}$$

$$\bar{I}_{i0} + x_i \leq Q_i \qquad \forall i \in \mathcal{I} \setminus \{|\mathcal{I}|\} \tag{130}$$

$$\sum_{i \in \mathcal{I} \setminus \{|\mathcal{I}|\}} x_i \leq \bar{I}_{|\mathcal{I}|0} \tag{131}$$

$$x_i \geq 0 \qquad \forall i \in \mathcal{I} \tag{132}$$

where the second-stage cost function $Q(x, \xi)$ is

$$Q(x, \xi) = \min \sum_{i \in \mathcal{I} \setminus \{|\mathcal{I}|\}} \left( t_{i,i+1} y_{i,i+1} + \frac{c_i}{Q_i} B_i^+ + c_i E_i^+ + p_i(-I_i^-) + p_i a_i + c_i b_i \right) \tag{133}$$

$$\text{s.t.} \quad y_{i,i+1} \leq C \qquad \forall i \in \mathcal{I} \setminus \{|\mathcal{I}|\} \tag{134}$$

$$I_{|\mathcal{I}|} = \bar{I}_{|\mathcal{I}|0} - \sum_{i \in \mathcal{I} \setminus \{|\mathcal{I}|\}} x_i + y_{|\mathcal{I}|-1,|\mathcal{I}|} \tag{135}$$

$$I_{|\mathcal{I}|} \leq \bar{I}_{|\mathcal{I}|0} \tag{136}$$

$$I_1 = \bar{I}_{|\mathcal{I}|0} + x_1 - \xi_1 - y_{1,2} \tag{137}$$

$$I_i = \bar{I}_{i0} + x_i - \xi_i + y_{i-1,i} - y_{i,i+1} \qquad \forall i \in \mathcal{I} \setminus \{1, |\mathcal{I}|\} \tag{138}$$

$$I_i^+ = \max\{0, I_i\} \qquad \forall i \in \mathcal{I} \setminus \{|\mathcal{I}|\} \tag{139}$$

$$I_i^- = \min\{0, I_i\} \qquad \forall i \in \mathcal{I} \setminus \{|\mathcal{I}|\} \tag{140}$$

$$E_i = I_i^+ - Q_i \qquad \forall i \in \mathcal{I} \setminus \{|\mathcal{I}|\} \tag{141}$$

$$E_i^+ = \max\{0, E_i\} \qquad \forall i \in \mathcal{I} \setminus \{|\mathcal{I}|\} \tag{142}$$

$$B_i = I_i^+ - x_i - \bar{I}_{i0} - E_i^+ \qquad \forall i \in \mathcal{I} \setminus \{|\mathcal{I}|\} \tag{143}$$

$$B_i^+ = \max\{0, B_i\} \qquad \forall i \in \mathcal{I} \setminus \{|\mathcal{I}|\} \tag{144}$$

$$I_i, B_i, E_i \in \mathbb{R} \qquad \forall i \in \mathcal{I} \setminus \{|\mathcal{I}|\} \tag{145}$$

$$y_{i,i+1}, I_i^+, B_i^+, E_i^+ \geq 0 \qquad \forall i \in \mathcal{I} \setminus \{|\mathcal{I}|\} \tag{146}$$

$$I_i^- \leq 0 \qquad \forall i \in \mathcal{I} \setminus \{|\mathcal{I}|\} \tag{147}$$

The formulated problem is similar to the model presented in Section 7.3 but with slight differences. The objective functions (127) and (133) aim to minimize the total expected costs, obtained by summing over all the penalties for delivery, re-balancing, extra and excess inventory, and stock-out in addition to the penalties of the slack variables. Constraints (128) encourage the target inventory quantity to be greater than or equal to the maximum number of consecutive withdrawals, while constraints (129) encourage the number of free racks to be greater than or equal to the maximum number of consecutive returns. To guarantee that a feasible solution can be found if the sum of the maximum number of consecutive withdrawn and returned bikes is greater than the station capacity the slack variables $a_i$ and $b_i$ allow for deviations from these quantities. Constraints (130) ensure that the sum of the allocated quantity and the initial availability at each station does not exceed the station's capacity. Constraints (131) indicate that the total number of allocated bikes within stations is less than the available quantity at the depot. Constraints (134) ensure that the number of bikes transported by the vehicle during re-balancing never exceeds its capacity. Constraints (135) guarantee that the quantity of bikes at the depot at the end of each scenario equals the initial bike availability plus the amount received from the last visited station, minus the quantities delivered to stations. Additionally, constraints (136) ensure that the number of bikes at the depot at the end of the re-balancing period does not exceed its capacity. The flow balance constraints for bikes at the first station on the route differ from the remaining stations. Specifically, constraints (137) ensure that for the first visited station, the amount of bikes at the end of the re-balancing period equals the sum of the initially available quantity and the quantity received from the depot, minus the amounts used to satisfy the demand and those bikes that are redistributed to subsequent stations on the route. Similarly, constraints (138) determine the inventory position at a station other than the first, as a function of the initial inventory level, the number allocated, the number withdrawn/returned, and the number redistributed to another station. Constraints (139) and (140) determine the surplus and stock-out quantities for each station. Constraints (141) and (142) calculate the number of bikes at each station that exceeds station capacity. Constraints (143) and (144) determine when more bikes are positioned at a station after re-balancing than were initially allocated, but not more than station capacity.

## E.1 Results

We now show the results of applying the methods in Section 6 for the previous two-stage problem. The uncertain parameter is the daily demand for bikes at each station, and our methodology uses additional covariate information to obtain a forecast of the demand. As discussed in Section 7, we estimate the parameters $g_i$ and $h_i$ in (128)-(129) as fractions of the demand in node $i$. In our experiments, we use the same features, data, and methodology used in Section 7.3 . To evaluate each method's optimality gap, we randomly split the dataset into training and test sets with an 80:20 ratio. We then calculate the out-of-sample cost by averaging the bike-sharing problem cost over the test data. This process is repeated with 10 sets of training and test data.

Figure 7 displays the out-of-sample costs for the proposed multiple replications procedure for the model discussed in Cavagnini et al. (2024). Despite similar variance in out-of-sample results across all methods, the approaches that consider multi-scenarios such as ER-SAA, KNN, and SAA outperform other methods that use point-wise forecast approximation, such as AD, CART, M5+AD, and LS, showing lower out-of-sample costs. Additionally, application-driven forecast methods within the point-wise forecast methods, such as AD and M5+AD, perform considerably better than methods such as CART and LS. Even though both problems (101)-(120) and (127)-(147) are relatively similar, the results differ due to the functional dependencies among the right-hand side presented because of constraints (128)-(129), as discussed in Section 3, indicating the need for a larger data set to solve this problem.

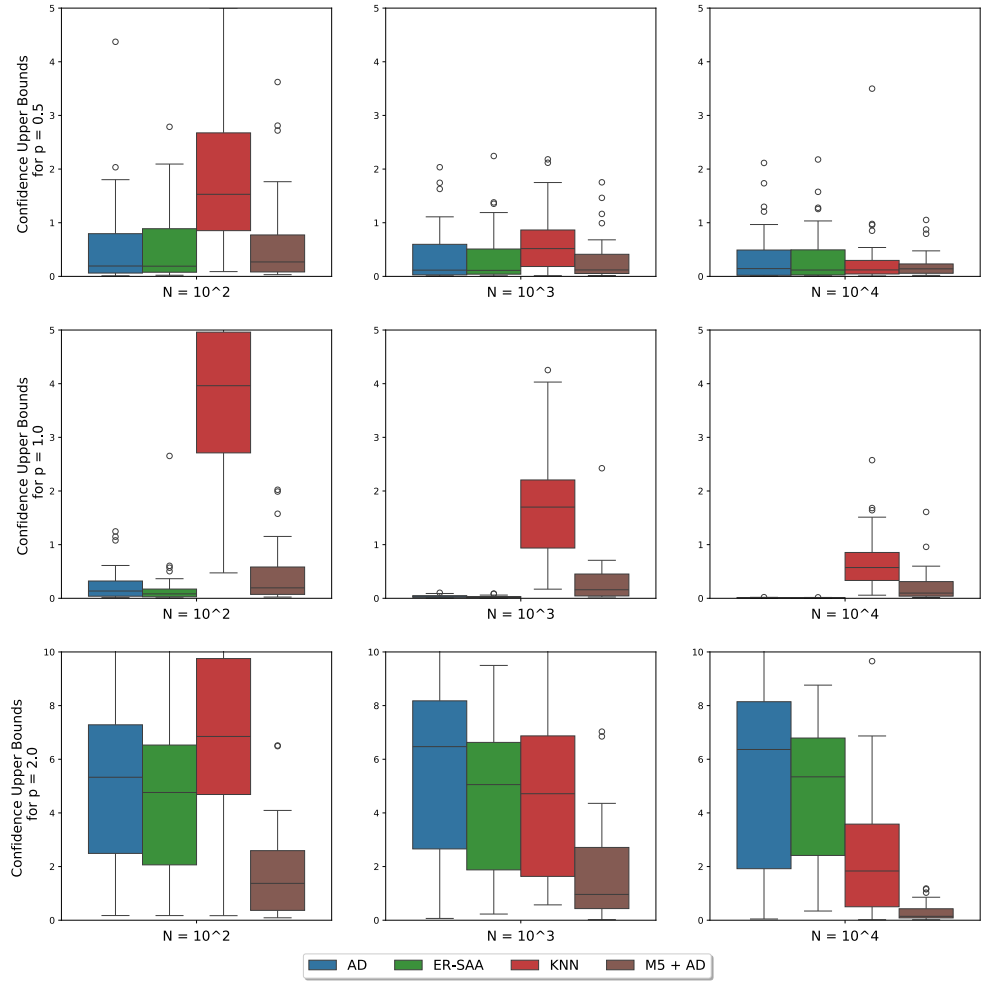

Figure 6: Detailed comparison of methods proposed in Section 6 for Problem 2 of shipment planning in Section 7.2.3. $p$ is the degree of the data generation procedure in (90), and $N$ is the number of samples with which each method is trained.

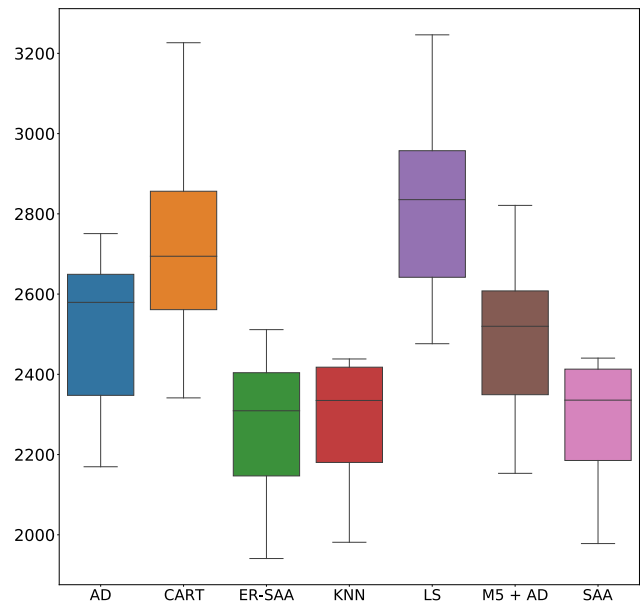

Figure 7: Average out-of-sample costs for 10 different train and test sets for methods proposed in Section 6.